\documentclass[runningheads]{svjour3}                   % onecolumn (standard format)
\smartqed  % flush right qed marks, e.g. at end of proof
\usepackage{graphicx}
\usepackage{amsmath}
\usepackage{epstopdf} % needed if you have eps figures in a pdflatex manuscript
\usepackage{mathptmx}      % use Times fonts if available on your TeX system
\usepackage{amssymb}
\usepackage[lined,ruled,norelsize]{algorithm2e}
\usepackage{bm}
\usepackage{graphicx}
\usepackage{color}
\usepackage[small,bf]{caption}
\usepackage{subfig}
\usepackage{hyperref}

\newcommand{\inv}[1]{{#1}^{-1}}

\newcommand{\pdesol}[1]{\underline{#1}}
\newcommand{\odesol}[1]{#1}
\newcommand{\disc}[1]{\uppercase{#1}}

\newcommand{\R}{{\mathbb{R}}}

\newcommand{\Qmat}{{\bm{Q}}}
\newcommand{\qmat}{{\bm{q}}}
\newcommand{\Smat}{{\bm{S}}}
\newcommand{\smat}{{\bm{s}}}
\newcommand{\Uvec}{{\bm{U}}}
\newcommand{\Fvec}{{\bm{F}}}

\newcommand{\tvec}{{\bm{t}}}
\newcommand{\rvec}{{\bm{r}}}

\newcommand{\tauvec}{{\bm{\tau}}}
\newcommand{\QLvec}[3]{\bm{#1}_{#2}^{#3}}
\newcommand{\QL}[4]{#1_{#2,#4}^{#3}}
\newcommand{\ULvec}[2]{\QLvec{U}{#1}{#2}}
\newcommand{\UL}[3]{\QL{U}{#1}{#2}{#3}}
\newcommand{\ULtmp}[3]{\QL{\tilde{U}}{#1}{#2}{#3}}
\newcommand{\FLvec}[2]{\QLvec{F}{#1}{#2}}
\newcommand{\FL}[3]{\QL{F}{#1}{#2}{#3}}
\newcommand{\BLvec}[2]{\QLvec{\tau}{#1}{}}

\makeatletter
\newcommand{\algorithmfootnote}[2][\footnotesize]{%
  \let\old@algocf@finish\@algocf@finish% Store algorithm finish macro
  \def\@algocf@finish{\old@algocf@finish% Update finish macro to insert "footnote"
    \leavevmode\rlap{\begin{minipage}{\linewidth}
    #1#2
    \end{minipage}}%
  }%
}
\makeatother

\renewcommand{\d}{\mathrm{d}}
\newcommand{\coloneq}{:=}

\journalname{BIT}
\begin{document}

\title{A multi-level spectral deferred correction method\thanks{Robert Speck and Daniel Ruprecht acknowledge supported by Swiss National Science Foundation grant 145271 under the lead agency agreement through the project "ExaSolvers" within the Priority Programme 1648 "Software for Exascale Computing" of the Deutsche Forschungsgemeinschaft. Matthias Bolten acknowledges support from DFG through the project "ExaStencils" within SPPEXA. Daniel Ruprecht and Matthew Emmett also thankfully acknowledge support by grant SNF-147597. Matthew Emmett and Michael Minion were  supported by the Applied Mathematics Program of the DOE Office of Advanced Scientific Computing Research under the U.S. Department of Energy under contract DE-AC02-05CH11231. Michael Minion was also supported by the U.S. National Science Foundation grant DMS-1217080.
}}


\author{Robert Speck \and Daniel Ruprecht \and Matthew Emmett \and Michael Minion \and Matthias Bolten \and Rolf Krause
}

\authorrunning{R. Speck, D. Ruprecht, M. Emmett, M. Minion, M. Bolten, R. Krause} % if too long for running head

\institute{R. Speck \at
	     J\"ulich Supercomputing Centre, Forschungszentrum J\"ulich, Germany and Institute of Computational Science, Universit{\`a} della Svizzera italiana, Lugano, Switzerland.\\
              \email{r.speck@fz-juelich.de}
            \and
           D. Ruprecht \at
	  Institute of Computational Science, Universit{\`a} della Svizzera italiana, Lugano, Switzerland.\\
	      \email{daniel.ruprecht@usi.ch}
              \and
              M. Emmett \at
	     Center for Computational Sciences and Engineering, Lawrence Berkeley National Laboratory, USA.\\
	      \email{mwemmett@lbl.gov}
              \and
              M. Minion \at
              Institute for Computational and Mathematical Engineering, Stanford University, USA.\\
              \email{mlminion@stanford.edu}
              \and
              M. Bolten \at
              Department of Mathematics, Bergische Universit\"at Wuppertal, Germany.\\
                \email{bolten@math.uni-wuppertal.de}
              \and
              R. Krause\at
              Institute of Computational Science, Universit{\`a} della Svizzera italiana, Lugano, Switzerland.\\
                \email{rolf.krause@usi.ch}
}

\date{Received: date / Accepted: date}

\maketitle

\begin{abstract}

The spectral deferred correction (SDC) method is an iterative scheme for computing a
higher-order collocation solution to an ODE by performing a series
of correction sweeps using a low-order timestepping method.
This paper examines a variation of SDC for the temporal integration of PDEs called
multi-level spectral deferred corrections (MLSDC), where sweeps are performed on a hierarchy of levels
and an FAS correction term, as in nonlinear multigrid methods, couples solutions on different levels.
Three different strategies to reduce the computational cost of correction
sweeps on the coarser levels are examined:
reducing the degrees of freedom, reducing the order of the spatial discretization,
and reducing the accuracy when solving linear systems arising in implicit temporal integration.
Several numerical examples demonstrate the effect of multi-level coarsening
on the convergence and cost of SDC integration.
In particular, MLSDC can provide significant savings in compute time compared to SDC for a three-dimensional problem.

%Include keywords and mathematical subject classification numbers as needed.
\end{abstract}

%% introductory section
\section{Introduction}
The numerical approximation of initial value ordinary differential equations is a fundamental problem in computational science, and many integration methods for problems of different character have been developed \cite{ascherPetzold,HairerI,HairerII}.
Among different solution strategies, this paper focuses on a class of iterative methods called Spectral Deferred Corrections (SDC) \cite{dutt2000spectral}, which is a variant of the defect and deferred correction methods developed in the 1960s \cite{bohmerHemkerStetter:1984,pereyra:1968,pereyra:1967,pereyra:1966,stetter:1974,zadunaisky:1964}.
In SDC methods, high-order temporal approximations are computed over a time\-step by discretizing and approximating a series of correction equations on intermediate substeps.
These corrections are applied iteratively to a provisional solution computed on the substeps, with each iteration -- or {\it sweep} -- improving the solution and raising the formal order of accuracy of the method, see e.g.~\cite{ChristliebEtAl2011_CMS,ChristliebEtAl2010_MoC,ShuEtAl2007}.
The correction equations are cast in the form of a Picard integral equation containing an explicitly calculated term corresponding to the temporal integration of the function values from the previous iteration.
Substeps in SDC methods are chosen to correspond to Gaussian quadrature nodes, and hence the integrals can be stably computed to a very high order of accuracy.

One attractive feature of SDC methods is that the numerical method used to approximate the correction equations can be low-order (even first-order) accurate, while the solution after many iterations can in principal be of arbitrarily high-order of accuracy.
This has been exploited to create SDC methods that allow the governing equations to be split into two or more pieces that can be treated either implicitly or explicitly and/or with different timesteps, see e.g.~\cite{bourlioux2003high,bouzarthMinion:2011,layton2004conservative,minion2003semi}.

For high-order SDC methods constructed from low-order propagators, the provisional solution and the solution after the first few correction iterations are of lower-order compared to the final solution.
Hence it is possible to reduce the computational work done on these early iterations by reducing the number of substeps (i.e. quadrature nodes) since higher-order integrals are not yet necessary.
In \cite{laytonMinion:2005,minion2003semi}, the number of substeps used in initial iterations of SDC methods is appropriately reduced to match the accuracy of the solution, and the methods there are referred to as {\it ladder methods}.
Ladder methods progress from a low-order coarse solution to a high-order fine solution by performing one or more SDC sweeps on the coarse level and then using an interpolated (in time and possibly space) version of the solution as the provisional solution for the next correction sweep.
In both \cite{laytonMinion:2005,minion2003semi} the authors conclude that the reduction in work obtained by using ladder methods is essentially offset by a corresponding decrease in accuracy, making ladder methods no more computationally efficient than non-ladder SDC methods.
On the other hand, in \cite{Layton2009}, SDC methods for a method of lines discretizations of PDEs are explored wherein the ladder strategy allows both spatial and temporal coarsening as well as the use of lower-order spatial discretizations in initial iterations.
The numerical results in \cite{Layton2009} indicate that adding spatial coarsening to SDC methods for PDEs can increase the overall efficiency of the timestepping scheme, although this evidence is based only on numerical experiments using simple test cases.

This paper significantly extends the idea of using spatial coarsening in SDC when solving PDEs.
A general multi-level strategy is analyzed wherein correction sweeps are applied to different levels as in the V-cycles of multigrid methods (e.g. \cite{brandt:1977,briggs}).
A similar strategy is used in the parallel full approximation scheme in space and time (PFASST), see~\cite{EmmettMinion2012,Minion2010} and also~\cite{SpeckEtAl2012}, to enable concurrency in time by iterating on multiple timesteps simultaneously.
As in nonlinear multigrid methods, multi-level SDC applies an FAS-type correction to enhance the accuracy of the solution on coarse levels.
Therefore, some of the fine sweeps required by a single-level SDC algorithm can be replaced by coarse sweeps, which are relatively cheaper when spatial coarsening strategies are used.
The paper introduces MLSDC and discusses three such spatial coarsening strategies: (1) reducing the number of degrees of freedom, (2) reducing the order of the discretization and (3) reducing the accuracy of implicit solves.
To enable the use of a high-order compact stencils for spatial operators, several modifications to SDC and MLSDC are presented that incorporate a weighting matrix.
It is shown for example problems in one and two dimensions that the number of MLSDC iterations required to converge to the collocation solution can be fewer than for SDC, even when the problem is poorly resolved in space.
Furthermore, results from a  three-dimensional benchmark problem demonstrate that MLSDC can significantly reduce time-to-solution compared to single-level SDC.

%% main section on MLSDC: SDC, MLSDC, PMG, compact stencils, first test
\section{Multi-level spectral deferred corrections} \label{sec:MLSDC}
The details of the MLSDC schemes are presented in this section.
The original SDC method is first reviewed in \S\ref{subsec:sdc}, while MLSDC along with a brief review of FAS corrections, the incorporation of weighting matrices and a discussion of different coarsening strategies is presented in \S\ref{subsec:mlsdc}.

\subsection{Spectral deferred corrections}\label{subsec:sdc}
SDC methods for ODEs were first introduced in \cite{dutt2000spectral}, and were subsequently refined and extended e.g. in~\cite{hansen2006convergence,huang2006accelerating,minion2003semi,minion2004semi}.
SDC methods iteratively compute the solution to the collocation equation by approximating a series of correction equations at spectral quadrature nodes using low-order substepping methods.
The derivation of SDC starts from the Picard integral form of a generic IVP given by
\begin{equation}
  \label{eq:picard}
  u(t) = u_{0} + \int_{0}^t f\bigl(u(s), s\bigr) \d s
\end{equation}
where $t \in [0,T]$, $u_0,u(t) \in \R^N$, and $f: \R^N \times \R \rightarrow \R^N$.
We now focus on a single timestep $[T_n, T_{n+1}]$, which is divided into substeps by defining
a set of quadrature nodes on the interval.  Here we consider Lobatto quadrature and denote
$M+1$ nodes  $\tvec := (t_m)_{m=0,\ldots,M}$ such that
$T_n=t_0 < t_{1} < \ldots < t_{M} = T_{n+1}$.
We now denote the collocation polynomial on $[T_{n}, T_{n+1}]$ by $u_p(t)$ and
write $U_j = u_p(t_j) \approx u(t_j)$.
In order to derive equations for the intermediate solutions $U_j$, we define quadrature weights
\begin{equation}
	\label{eq:quad_weights}
	q_{m,j} \coloneq \frac{1}{\Delta t} \int_{T_{n}}^{t_{m}} l_{j}(s) \ ds, \ m=0,\ldots,M, \ j=0,\ldots,M
\end{equation}
where $(l_{j})_{j=0,\ldots,M}$ are the Lagrange polynomials defined by the nodes $\tvec$, and $\Delta t = T_{N+1}-T_N$.
Inserting $u_p(t)$ into~\eqref{eq:picard} and noting that the quadrature with weights defined in~\eqref{eq:quad_weights} integrates the polynomial $u_p(t)$ exactly, we obtain
\begin{equation}
	\label{eq:disc_coll}
	U_m = u_0 + \Delta t \sum_{j=0}^{M} q_{m,j} f(U_{j}, t_j), \ m=0, \ldots, M.
\end{equation}
For a more compact notation, we now define the {\it integration matrix} $\qmat$ to be the $M+1 \times M+1$
matrix consisting of entries $q_{m,j}$.  Note that because we use Gauss-Lobatto nodes, the first row of $\qmat$ is all zeros.  Next, we  denote
\begin{equation}
	\Uvec := \left[ U_0, \ldots, U_M \right]^T,\nonumber
\end{equation}
and
\begin{equation}
 	\Fvec(\Uvec) := \left[ F_0, \ldots, F_M \right]^T := \left[ f(U_0, t_0), \ldots, f(U_M, t_M) \right]^T. \nonumber
 \end{equation}
In order to multiply the integration matrix $\qmat$ with the vector of the right-hand side values, we define $\Qmat \coloneq \qmat\otimes \bm{I}_N$ where $\bm{I}_N\in\mathbb{R}^{N\times N}$ is the identity matrix and $\otimes$ is the Kronecker product.
With these definitions, the set of equations in~\eqref{eq:disc_coll} can be written more compactly as
\begin{equation}
  \label{eq:compact}
  \Uvec = \Uvec_0 + \Delta t\, \Qmat\, \Fvec(\Uvec) \nonumber
\end{equation}
where $\Uvec_0 := U_0 \otimes \bm{I}_N$.
Eq. \eqref{eq:compact} is an implicit equation for the unknowns in $\Uvec$, and is also referred to
as the collocation formulation.
Because we use Gauss-Lobatto nodes, the value $U_M$ readily approximates the solution $u(T_{n+1})$.

Here, we consider ODEs that can be split into stiff ($f^I$) and non-stiff ($f^E$) pieces so that
\begin{equation}
 f(u(t),t) = f^E\bigl(u(t), t\bigr) + f^I\bigl(u(t), t\bigr).\nonumber
\end{equation}
SDC iterations begin by spreading the initial condition $U_0$ to each of the collocation nodes so that the provisional solution $\Uvec^0$ is given by $\Uvec^0 = [U_0, \cdots, U_0]$.
We define by
\begin{equation}
	s_{m,j} \coloneq \frac{1}{\Delta t} \int_{t_{m-1}}^{t_m} l_{j}(s) \ ds, \ m=1, \ldots, M\nonumber
\end{equation}
the quadrature weights for node-to-note integration, approximating integrals over $[t_{m-1}, t_{m}]$, and as $\smat$ the $M \times M+1$ matrix consisting of the entries $s_{m,j}$.
Note that $\smat$ can be easily constructed from the integration matrix $\qmat$.
Furthermore, we denote as before $\Smat\coloneq\smat\otimes\bm{I}_N$.
Then, the semi-implicit update equation corresponding to the forward/backward Euler substepping method for computing $\Uvec^{k+1}$ is given by
\begin{multline}
  \label{eq:imexsdc}
  U^{k+1}_{m+1} = U^{k+1}_m
    + \Delta t_m
      \bigl[ f^E(U^{k+1}_{m}, t_{m}) - f^E(U^k_{m}, t_{m}) \bigr] \\
    + \Delta t_m
      \bigl[ f^I(U^{k+1}_{m+1}, t_{m+1}) - f^I(U^k_{m+1}, t_{m+1}) \bigr]
    + \Delta t\, S^{k}_{m}
\end{multline}
where $S^{k}_{m}$ is the $m^{\rm th}$ row of $\Smat \Fvec(\Uvec^k)$ and $\Delta t_{m} \coloneq t_{m+1} - t_{m}$.
The process of solving \eqref{eq:imexsdc} at each node is referred to as an \emph{SDC sweep} or an \emph{SDC iteration} (see Algorithm~\ref{alg:sdcsweep}).
SDC with a fixed number of $k$ iterations and first-order sweeps is formally $O(\Delta t^k)$ up to the accuracy
of the underlying integration rule~\cite{ChristliebEtAl2009,ShuEtAl2007}.
When SDC iterations converge, the scheme becomes equivalent to the collocation scheme determined by the quadrature nodes, and hence is of order $2M$ with $M+1$ Lobatto nodes.
\begin{algorithm}[t]
  \algorithmfootnote{The FAS correction, denoted by $\tauvec$, is included here to ellucidate how FAS corrections derived in \S\ref{subsec:mlsdc} are incorporated into an SDC sweep -- for plain, single level SDC algorithms the FAS correction $\tauvec$ would be zero.}
  \SetKwComment{Comment}{\# }{}
  \SetCommentSty{textit}
  \DontPrintSemicolon

  \KwData{Initial $U_0$, function evaluations $\Fvec(\Uvec^k)$ from the previous iteration, and (optionally) FAS corrections $\tauvec$.}
  \KwResult{Solution $\Uvec^{k+1}$ and function evaluations $\Fvec(\Uvec^{k+1})$.}

  \BlankLine
  \Comment{Compute integrals}
  \For{$m=0 \ldots M-1$}{
    $S^{k}_{m} \longleftarrow \Delta t \sum_{j=0}^M s_{m,j} (F^{E,k}_{j} + F^{I,k}_j)$
  }

  \BlankLine
  \Comment{Set initial condition and compute function evaluation}
  $t \longleftarrow t_0$; $U^{k+1}_0 \longleftarrow U_0$ \;
  $F^{E,k+1}_0 \longleftarrow f^E(U_0, t)$ \;
  $F^{I,k+1}_0 \longleftarrow f^I(U_0, t)$ \;

  \BlankLine
  \Comment{Forward/backward Euler substepping for correction}
  \For{$m=0 \ldots M-1$}{
    $t \longleftarrow t + \Delta t_m$ \;
    ${\rm RHS} \longleftarrow U^{k+1}_{m} + \Delta t_m \bigl( F^{E,k+1}_{m} - F^{E,k}_{m} - F^{I,k}_{m+1} \bigr) + S^k_{m} + \tau_m$ \;
    $U^{k+1}_{m+1} \longleftarrow {\rm Solve}\bigl( U - \Delta t_m f^I(U, t) = {\rm RHS} \bigr)$ for $U$ \;
    $F^{E,k+1}_{m+1} \longleftarrow f^E(U^{k+1}_{m+1}, t)$ \;
    $F^{I,k+1}_{m+1} \longleftarrow f^I(U^{k+1}_{m+1}, t)$ \;
  }

  \caption{IMEX SDC sweep algorithm.}
  \label{alg:sdcsweep}
\end{algorithm}

%\subsubsection{SDC for stiff problems}
It has been shown \cite{huang2006accelerating,laytonMinion:2005}
that in certain situations (particularly
stiff equations)  the convergence of SDC iterates can
slow down considerably for large values of $\Delta t$.  For a fixed
number of iterations, this lack of convergence is characterized
by order reduction.  Hence in this study, to allow for a reasonable comparison of SDC and MLSDC, we perform iterations until a specified convergence criterion is met.
Convergence is monitored  by computing the SDC residual
\begin{equation}
  \label{eq:residual}
  \rvec^k = \Uvec_0 + \Delta t \Qmat \Fvec(\Uvec^k) - \Uvec^k,
\end{equation}
and the iteration is terminated when the norm of the residual drops
below a prescribed tolerance.
Similary, if SDC or MLSDC are used to solve the collocation problem up to some fixed tolerance, one also observes a significant increase in the number of iterations required to reach a set tolerance.
Accelerating the convergence of SDC for stiff problems has been studied in e.g.~\cite{HuangEtAl2006,Weiser2013}.

%% subsection on MLSDC
\subsection{Multi-level spectral deferred corrections}\label{subsec:mlsdc}
In multi-level SDC (MLSDC), SDC sweeps are performed on a hierarchy of discretizations
or \emph{levels}  to solve the collocation equation~\eqref{eq:compact}.
This section presents the details of the MLSDC iterations
for a generic set of levels, and in Sect. \ref{sec:mlsdc_spatial_coarsening},
three different coarsening strategies are explored.
For the following, we define levels $\ell=1 \ldots L$, where $\ell = 1$ is
the discretization that is to be solved (referred to generically
as the {\it fine} level), and subsequent
levels $\ell=2 \ldots L$ are defined by successive coarsening
of a type to be specified later.

\subsubsection{FAS correction}
Solutions on different MLSDC levels are coupled in the same manner as used in the full approximation scheme (FAS) for nonlinear multigrid methods (see e.g. \cite{brandt:1977}).
The FAS correction for coarse SDC iterations
is determined by considering SDC as an iterative method for solving
the collocation formulation~\eqref{eq:compact}, where the
operators $A_\ell$ are given by $A_\ell(\Uvec_\ell) \equiv
\Uvec_\ell - \Delta t \Qmat_\ell \Fvec_\ell(\Uvec_\ell)$.
Note that the approximations $A_\ell$ of the operator $A$ can differ substantially between levels as will be discussed in \S\ref{sec:mlsdc_spatial_coarsening}.
Furthermore, we assume that suitable restriction (denote by $R$) and interpolation operators between levels are available, see \S\ref{subsubsec:transfer}.
The FAS correction for coarse-grid sweeps is  then given by
\begin{equation}
  \label{eq:tau_sdc}
  \bm{\tau}_{\ell+1} = A_{\ell+1}(R \Uvec_{\ell}) - R A_\ell(\Uvec_\ell) =  \Delta t \bigl( R \Qmat_{\ell} \Fvec_{\ell}(\Uvec_\ell)
                           - \Qmat_{\ell+1} \Fvec_{\ell+1} (R \Uvec_\ell)\bigr).
\end{equation}
In particular, if the fine residual is zero (i.e., $\Uvec_{\ell} \equiv
\Uvec_{0,\ell} + \Delta t \Qmat_{\ell} \Fvec_{\ell}(\Uvec_\ell)$) the FAS-corrected
coarse equation becomes
\begin{eqnarray}
  \Uvec_{\ell+1}- \Delta t \Qmat_{\ell+1} \Fvec_{\ell+1}(\Uvec_{\ell+1})
   &  = & R\Uvec_{0,\ell} + \Delta t \bigl( R \Qmat_{\ell} \Fvec_{\ell}(\Uvec_\ell)
                                - \Qmat_{\ell+1} \Fvec_{\ell+1}(R\Uvec_\ell) \bigr) \nonumber\\
  &   =  & R \Uvec_{\ell} - \Delta t  \Qmat_{\ell+1} \Fvec_{\ell+1}(R\Uvec_\ell)\nonumber
\end{eqnarray}
so that the coarse solution is the restriction of the fine solution. Note that for multi-level schemes, FAS-corrections from finer levels need to be restricted and incorporated to coarser levels as well, i.e.~if on level $\ell$ the equation is already corrected by $\bm{\tau}_\ell$ with
\begin{equation}
  A_\ell(\Uvec_\ell) = \Uvec_\ell - \Delta t \Qmat_\ell \Fvec_\ell(\Uvec_\ell) - \bm{\tau}_\ell,\nonumber
\end{equation}
the correction $\bm{\tau}_{\ell+1}$ for level $\ell+1$ is then given by
\begin{equation}
  \bm{\tau}_{\ell+1} = A_{\ell+1}(R \Uvec_{\ell}) - R A_\ell(\Uvec_\ell) =  \Delta t \bigl( R \Qmat_{\ell} \Fvec_{\ell}(\Uvec_\ell)
                           - \Qmat_{\ell+1} \Fvec_{\ell+1} (R \Uvec_\ell)\bigr) + R\bm{\tau}_\ell.\nonumber
\end{equation}
Coarse levels thus include the FAS corrections of all finer levels.

\subsubsection{The MLSDC algorithm}
The MLSDC scheme introduced here proceeds as follows.  The initial
condition $U_0$ and its function evaluation are spread to each of the
collocation nodes on the finest level so that the first provisional
solution $\Uvec^0_1$ is given by
\begin{equation}
\Uvec^0_1 = [ U_0, \ldots, U_0 ].\nonumber
\end{equation}
A single MLSDC iteration then consists of the following steps:
\begin{enumerate}
\item Perform one fine SDC sweep using the values $\Uvec^{k}_1$ and
  $\Fvec_1(\Uvec^{k}_1)$.  This will yield provisional updated values
  $\Uvec^{k+1}_1$ and $\Fvec_1(\Uvec^{k+1}_1)$.
\item Sweep from fine to coarse: for each $\ell=2\ldots L$:
  \begin{enumerate}
  \item Restrict the fine values $\Uvec^{k+1}_{\ell-1}$ to the coarse
    values $\Uvec_{\ell}^{k}$ and compute $\Fvec_{\ell}(\Uvec_{\ell}^{k})$.
  \item Compute the FAS correction
    $\tauvec^k_{\ell}$ using $\Fvec_{\ell-1}(\Uvec^{k+1}_{\ell-1})$, $\Fvec_{\ell}(\Uvec^{k}_\ell)$,
    and $\tauvec^k_{\ell-1}$ (if available).
  \item Perform $n_{\ell}$ SDC sweeps with the values on level $\ell$
    beginning with $\ULvec{\ell}{k}$, $\Fvec_{\ell}(\ULvec{\ell}{k})$ and the FAS correction
    $\tauvec^k_{\ell}$.  This will yield new values $\ULvec{\ell}{k+1}$
    and $\Fvec_{\ell}(\ULvec{\ell}{k+1})$.
  \end{enumerate}
\item Sweep from coarse to fine: for each $\ell=L-1\ldots 1$:
  \begin{enumerate}
  \item Interpolate coarse grid correction $\ULvec{\ell+1}{k+1} - R \ULvec{\ell}{k+1}$ and
    add to $\ULvec{\ell}{k+1}$.  Recompute new values
    $\Fvec_{\ell}(\ULvec{\ell}{k+1})$
  \item If $\ell > 1$, perform $n_{\ell}$ SDC sweeps beginning with
    values $\ULvec{\ell}{k+1}$, $\Fvec_{\ell}(\ULvec{\ell}{k+1})$ and the FAS
    correction $\tauvec^k_{\ell}$.  This will once again yield new
    values $\ULvec{\ell}{k+1}$ and $\Fvec_{\ell}(\ULvec{\ell}{k+1})$.
  \end{enumerate}
\end{enumerate}
Note that when interpolating from coarse to fine levels the correction
$\Uvec^{k+1}_{\ell+1} - R \Uvec^k_{\ell+1}$ is interpolated and
subsequently added to $\Uvec^{k+1}_{\ell}$ instead of simply
overwriting the fine values with interpolated coarse values.  Also
note that instead of interpolating solution values
$\Uvec^{k+1}_{\ell+1}$ to $\Uvec^{k+1}_\ell$ and immediately
re-evaluating the function values $\Fvec_\ell(\Uvec^{k+1}_\ell)$, the
change in the function values can be interpolated as well.
Doing so reduces the cost of
the interpolation step, but possibly at the cost of
increasing the number of MLSDC iterations required to reach
convergence. Since no significant increase could be observed during our tests,
we skip the re-evaluation of the right-hand side and use interpolation of the
coarse function values throughout this work.
The above is summarized by Algorithm~\ref{alg:mlsdc}.

\begin{algorithm}[t]
  \SetKwComment{Comment}{\# }{}
  \SetCommentSty{textit}
  \DontPrintSemicolon

  \KwData{Initial $\UL{1}{k}{0}$ and function evaluations $\FLvec{1}{k}$ from the previous iteration on the fine level.}
  \KwResult{Solution $\ULvec{\ell}{k+1}$ and function evaluations $\FLvec{\ell}{k+1}$ on all levels.}

  \BlankLine
  \Comment{Perform fine sweep and check convergence criteria}
  $\ULvec{1}{k+1}$, $\FLvec{1}{k+1} \longleftarrow$ SDCSweep$\bigl(\ULvec{1}{k},\,\FLvec{1}{k}\bigr)$ \;
  \If{ fine level has converged }{
    return \;
  }

  \BlankLine
  \Comment{Cycle from fine to coarse}
  \For{$\ell=1 \ldots L-1$}{
    \Comment{Restrict, re-evaluate, and save restriction (used later during interpolation)}
    \For{$m = 0 \ldots M$}{
      $\UL{\ell+1}{k}{m} \longleftarrow$ Restrict$\bigl( \UL{\ell}{k+1}{m} \bigr)$ \;
      $\FL{\ell+1}{k}{m} \longleftarrow$ FEval$\bigl(\UL{\ell+1}{k+1}{m} \bigr)$ \;
      $\ULtmp{\ell+1}{k}{m} \longleftarrow \UL{\ell+1}{k}{m}$\;
    }
    \Comment{Compute FAS correction and sweep}
    $\BLvec{\ell+1}{k} \longleftarrow$ FAS$\bigl(\FLvec{\ell}{k+1},\, \FLvec{\ell+1}{k},\, \BLvec{\ell}{k}\bigr)$ \;
    $\ULvec{\ell+1}{k+1}$, $\FLvec{\ell+1}{k+1} \longleftarrow$ SDCSweep$\bigl(\ULvec{\ell+1}{k},\,\FLvec{\ell+1}{k},\, \BLvec{\ell+1}{k} \bigr)$ \;
  }

  \BlankLine
  \Comment{Cycle from coarse to fine}
  \For{$\ell=L-1 \ldots 2$}{
    \Comment{Interpolate coarse correction and re-evaluate}
    \For{$m = 0 \ldots M$}{
      $\UL{\ell}{k+1}{m} \longleftarrow \UL{\ell}{k+1}{m} +$ Interpolate$\bigl(\UL{\ell+1}{k+1}{m} - \ULtmp{\ell+1}{k}{m} \bigr)$ \;
      $\FL{\ell}{k+1}{m} \longleftarrow$ FEval$\bigl(\UL{\ell}{k+1}{m}\bigr)$\;
    }
    $\ULvec{\ell}{k+1}$, $\FLvec{\ell}{k+1} \longleftarrow$ SDCSweep$\bigl(\ULvec{\ell}{k+1},\,\FLvec{\ell}{k+1},\, \BLvec{\ell}{k} \bigr)$ \;
  }

  \BlankLine
  \Comment{Return to finest level before next iteration}
    \For{$m = 0 \ldots M$}{
      $\UL{1}{k+1}{m} \longleftarrow \UL{1}{k+1}{m} +$ Interpolate$\bigl(\UL{2}{k+1}{m} - \ULtmp{2}{k}{m} \bigr)$ \;
      $\FL{1}{k+1}{m} \longleftarrow$ FEval$\bigl(\UL{1}{k+1}{m}\bigr)$\;
    }

  \BlankLine
  \caption{MLSDC iteration for $L$ levels.}
  \label{alg:mlsdc}
\end{algorithm}

\subsubsection{Semi-implicit MLSDC with compact stencils}\label{subsec:imex_linear}
In order to achieve higher-order accuracy with finite difference discretizations in space, the use of Mehr\-stellen discretizations is a common technique especially when using multigrid methods~\cite{trottenberg_multigrid:_2000}. While the straightforward use of larger stencils leads to larger matrix bandwidths and higher communication costs during parallel runs, \emph{high-order compact} schemes allow for high-order accuracy with stencils of minimal extent~\cite{spotz_high-order_1996}. The compact stencil for a given discretization is obtained by approximating the leading order error term by a finite difference approximation of the right-hand side, resulting in a weighting matrix.
Discretizing e.g.~the heat equation $\pdesol{u}_t = \nabla^2 \pdesol{u}$ in space\footnote{We adopt here and in the upcoming examples the following notation: Solutions of PDEs are
denoted with an underline, e.g.  $\pdesol{u}$, and depend continuously on
one or more spatial variables and a time variable.  Discretizing a
PDE in space by the method of lines results in an IVP with dimension $N$ equal to the
degrees of freedom of the spatial discretization.  The solution of
such an IVP is a vector-valued function denoted by a lower case letter, e.g.~$\odesol{u}$, and
depends continuously on time. The numerical approximation of $\odesol{u}$
at some point in time $t_{m}$ is
denoted by a capital letter, e.g.~$\disc{U}_{m}^{k}$, where $k$ corresponds to the iteration number.}
yields
\begin{align}
	Wu_t = Au\nonumber
\end{align}
with system matrix $A$ and weighting matrix $W$. Formally, the discrete Laplacian is given by $\inv{W}A$. Using this approach, a fourth-order approximation of the Laplacian can be achieved using only nearest neighbors (three-point stencil in 1D, nine-point-stencil in 2D, 19-point stencil in 3D).
For further reading on compact schemes we refer to~\cite{lele_compact_1992,spotz_high-order_1996,trottenberg_multigrid:_2000}.

The presence of a weighting matrix requires some modifications to MLSDC.  We start with the semi-implicit SDC update equation~\eqref{eq:imexsdc} given by
\begin{multline}
\label{eq:imexsdc_ii}
  U^{k+1}_{m+1} = U^{k+1}_m
    + \Delta t_m
      \bigl[ f^E(U^{k+1}_{m}, t_{m}) - f^E(U^k_{m}, t_{m}) \bigr] \\
    + \Delta t_m
      \bigl[ f^I(U^{k+1}_{m+1}, t_{m}) - f^I(U^k_{m+1}, t_{m}) \bigr]
    + \Delta t\, S^{k}_m.
\end{multline}
Next, we assume a linear, autonomous implicit part $f^I(U,t) = f^I(U)= \inv{W}AU$ for a spatial vector $U$ with sparse matrices $W$ and $A$ stemming from the discretization of the Laplacian with compact stencils. Furthermore, we define
\begin{equation}
	\tilde{f}^I(U) = AU\nonumber
\end{equation}
so that
\begin{equation}
\label{eq:f_tilde_def}
 \tilde{f}^I(U) = Wf^I(U).
\end{equation}
With these definitions \eqref{eq:imexsdc_ii} becomes
\begin{multline}
\left(I-\Delta t_m\, \inv{W}A\right)U^{k+1}_{m+1} = U^{k+1}_m
    + \Delta t_m
      \bigl[ f^E(U^{k+1}_{m}, t_{m}) - f^E(U^k_{m}, t_{m}) \bigr] \nonumber\\ - \Delta t_m\, \inv{W}AU^k_{m+1} + \Delta t\, S^{k}_m.\nonumber
\end{multline}
Since the operator $\left(I-\Delta t_m\, \inv{W}A\right)$ is not sparse, we
avoid computing with it by multiplying the equation above by $W$, so that
\begin{multline}
\label{eq:linear_imex_sdc_W}
\left(W-\Delta t_m\, A\right)U^{k+1}_{m+1} = WU^{k+1}_m
    + \Delta t_m
      W\bigl[ f^E(U^{k+1}_{m}, t_{m}) - f^E(U^k_{m}, t_{m}) \bigr] \\ - \Delta t_m\, \tilde{f}^I(U^k_{m+1}) + \Delta t\, \tilde{S}^{k}_m
\end{multline}
where $\tilde{S}^{k}_m$ now represents the $m^{\rm th}$ row of $\Smat \tilde{\Fvec}^k(\Uvec^k)$, using $Wf^E(U^k_m,t_{m})$ and $\tilde{f^I}(U^k_m)$ instead of $f^E(U^k_m,t_{m})$ and $f^I(U^k_m)$ as integrands, that is $\tilde{S}^{k}_m = \sum_{j=0}^M s_{m,j} \bigl( W f^E(U^k_j,t_{j}) + \tilde{f^I}(U^k_j) \bigr)$.

While this equation avoids the inversion of $W$, the computation of the residual does not. By equation~\eqref{eq:residual}, the $m^{\rm th}$ component of the residual at iteration $k$ reads either
\begin{equation}
  r^k_m = U_0 + \Delta t  \left(\Qmat \Fvec(\Uvec^k) \right)_m- U^k_m,\nonumber
\end{equation}
or, after multiplication with $W$,
\begin{equation}
  Wr^k_m = WU_0 + \Delta t \left( \Qmat \tilde{\Fvec}(\Uvec^{k})\right)_m - WU_m^k.\nonumber
\end{equation}
Both equations require the solution of a linear system with matrix $W$, either to compute the components of $\Fvec(\Uvec^k)$ from~\eqref{eq:f_tilde_def} or to retrieve $r^k_m$ from $W r^k_m$.
Note that the subscript $m$ denotes here the $m^{\rm th}$ column.
Thus, we either need to obtain $r^k_m$ from $Wr^k_m$ (in case $Wf^E$ is stored during the SDC sweep) or $f^I$ from $\tilde{f}^I$ (in case $f^E$ is stored). In either case, solving a linear system with the weighting matrix becomes inevitable for the computation of the formally correct residual.

Furthermore, evaluating \eqref{eq:tau_sdc} for the FAS correction also requires the explicit use of $f^E$ and $f^I = \inv{W}\tilde{f}^I$ to compute $R \Qmat_{\ell} \Fvec_{\ell}(\Uvec_{\ell})$.
Moreover, from \eqref{eq:linear_imex_sdc_W} we note that weighted SDC sweeps on coarse levels $\ell+1$ require the computation of $W_{\ell+1}\tau_{\ell+1,m}$ on all coarse nodes ${\bm t}_\ell$ so that $\Qmat_{\ell+1} \Fvec_{\ell+1}(R\Uvec_{\ell})$ can be replaced by $\Qmat_{\ell+1} \tilde{\Fvec}_{\ell+1}(R \Uvec_{\ell})$.
For spatial discretizations in which both parts $f^E$ and $f^I$ of the right-hand side make use of weighting matrices $W^E$ and $W^I$ or e.g.~for finite element discretizations with a mass matrix, we note that similar modifications to the MLSDC scheme as presented here must be made.
The investigation of MLSDC for finite element discretizations is left for future work.

\subsubsection{Coarsening strategies}\label{sec:mlsdc_spatial_coarsening}
The goal in MLSDC methods is to reduce the total cost of the method by performing
SDC sweeps on coarsened levels at reduced computational cost. In this section
we describe the three types of spatial coarsening used in the numerical examples:\vspace{0.5\baselineskip}
\begin{enumerate}
	\itemsep1em
	\item {\sc Reduced resolution in space}: Use fewer degrees of freedom for the spatial representation (e.g. nodes, cells, points, particles, etc.) on the coarse levels.
This directly translates into significant computational savings for evaluations of $f$, particularly for 3D problems.
This approach requires spatial interpolation and restriction operators to transfer the solution between levels.
	\item {\sc Reduced order in space}: Use a spatial discretization on the coarse levels that is of reduced order.
Lower-order finite difference stencils, for example, are typically cheaper to evaluate than higher-order ones, see~\cite{RuprechtKrause2012} for an application of this strategy for the time-parallel Parareal method.
	\item {\sc Reduced implicit solve in space:} Use only a few iterations of a spatial solver in every substep, if an implicit or implicit-explicit method is used in the SDC sweeps.
	By not solving the linear or nonlinear system in each SDC substep to full accuracy, savings in execution time can be achieved.
\end{enumerate}
\vspace{0.5\baselineskip}

We note that a fourth possibility not pursued here is to use a
simplified physical representation of the problem on coarse levels.
This approach requires a detailed understanding of the
problem to derive suitable coarse level models and appropriate coarsening
and interpolation operators.
Similar ideas have been studied for Parareal in~\cite{DaiEtAl2013_ESAIM,HautWingate2013}.

The spatial coarsening strategies outlined above can significantly
reduce the cost of a coarse level SDC substep, but do not affect the
number of substeps used.  In principle, it is also possible to
reduce the number of quadrature nodes on coarser levels as in the ladder
schemes mentioned in the introduction.
In this paper, no such temporal coarsening is applied and we
focus on the application of spatial coarsening strategies which leads
to a large reduction of the runtime for coarse level sweeps.

\subsubsection{Transfer operators}\label{subsubsec:transfer}
In order to apply Strategy 1 and reduce the number of spatial degrees of freedom, transfer operators between different levels are required.
In the tests presented here that are based on finite difference discretizations on simple cartesian meshes, the spatial degrees of freedom are aligned, so that simple injection can be used for restriction.

We have observed that the order of the used spatial interpolation has a strong impact on the convergence of MLSDC.
While global information transfer when using e.g.~spectral methods does not influence the convergence properties of MLSDC, the use of local Lagrangian interpolation for finite difference stencils has to be applied with care.
In numerical experiments not documented here, MLSDC with simple linear interpolation required twice as many iterations as MLSDC with fifth-order spatial interpolation.
Further, low resolutions in space combined with low-order interpolation led to significant degradation of the convergence speed of MLSDC, while high spatial resolutions were much less sensitive.
Throughout the paper, Strategy 1 is applied with third-order Lagrangian interpolation, which has proven to be sufficient in all cases studied here.

We note that the transfer operators would be different if e.g.~finite elements were used and operators between element spaces of different order and/or on different meshes would be required.

\subsubsection{Stability of SDC and MLSDC}
Stability domains for SDC are presented in e.g.~\cite{dutt2000spectral}.
The stability of semi-implicit SDC is addressed in~\cite{minion2003semi} and the issue of order reduction for stiff problems is discussed.
Split SDC methods are further analyzed theoretically and numerically in~\cite{HagstromZhou2006}.
A stability analysis for MLSDC is complicated by the fact that it would need to consider the effects of the different spatial coarsening strategies laid out in~\ref{sec:mlsdc_spatial_coarsening}.
Therefore, it cannot simply use Dahlquist's test equation but has to resort to some well-defined PDE examples in order to assess stability.
Hence, for MLSDC the results presented here are experimental but development of a theory for the convergence properties of MLSDC is ongoing work.
However, in all examples presented below, stability properties of SDC and MLSDC appeared to be comparable, but a comprehensive analysis is left for future work.

%% section with examples
\section{Numerical Examples}\label{sec:num_examples}
In this section we investigate the performance of MLSDC for four numerical examples.
%Since MLSDC requires only one fine SDC sweep per MLSDC iteration, if the number of MLSDC
%iterations required to converge to a given tolerance is less than the corresponding number of
%fine SDC iterations, then
%an overall savings in computational cost can be achieved if the computational cost
%of the coarse levels is sufficiently small.  e.
First, in order to demonstrate that the FAS correction in MLSDC is not unusable for hyperbolic problems per se, performance for the 1D wave equation is studied in \S\ref{subsec:wave_eq}.
To investigate performance for a nonlinear problem, MLSDC is then applied to the 1D viscous Burgers' equation in \S\ref{subsec:visc_burg}.
A detailed investigation of different error components is given and we verify that the FAS corrections allow the solutions on coarse levels to converge to the accuracy determined by the discretization on the \emph{finest} level.
The 2D Navier-Stokes equations in vorticity-velocity form are solved in \S\ref{subsec:shear_layer}, showing again a reduction of the number of required iterations by MLSDC, although using a coarsened spatial resolution is found to have a negative impact on convergence, if the fine level is already under-resolved.
In \S\ref{subsec:visc_burger3d}, a {\sc fortran} implementation of MLSDC is applied to the three-dimensional Burgers' equation and it is demonstrated that the reduction in fine level sweeps translates into a significant reduction of computing time.
Throughout all examples, we make use of a linear geometric multigrid solver~\cite{chow2006,trottenberg_multigrid:_2000} with  JOR relaxation in 3D and SOR relaxation 1D and 2D as smoothers, to solve the linear problems in the implicit part as well as to solve the linear system with the weighting matrix for the residual and the FAS correction.
%As relaxation parameter on each level the maximum absolute column sum of the respective system matrix is chosen. This guarantees convergence according to \cite{art:ARIC07}.
The parallel implementation of the multigrid solver used for the last example is described in~\cite{bolten:2014}.
%\todo{Matthias: is this correct? any information missing?}

In the examples below, we compare the number of sweeps on the fine and most expensive level required by SDC or MLSDC to converge up to a set tolerance.
For SDC, which sweeps only on the fine level, this number is identical to the number of iterations.
For MLSDC, each iteration consists of one cycle through the level hierarchy, starting from the finest level, going up to the coarsest and then down again, with one SDC sweep on each level on the way up and down, cf. Algorithm~\ref{alg:mlsdc}.
Except for the last iteration, the final fine sweep is also the first fine sweep of the next iteration, so that for MLSDC the number of fine sweeps is equal to the number of iterations plus one.
Note that a factor of two coarsening in the spatial resolution in each dimension yields a factor of eight reduction in degrees of freedom in three dimensions, which makes coarse level sweeps significantly less expensive.

\subsection{Wave equation}\label{subsec:wave_eq}
For spatial multigrid, the FAS formalism is mostly derived and analyzed for stationary elliptic or parabolic problems, although there are examples of applications to hyperbolic problems as well~\cite{Alam2006,south1977}.
Here, as a first test, we investigate the performance of MLSDC for a simple 1D wave equation to verify that the FAS procedure as used in MLSDC does not break down for a hyperbolic problem per se.
The problem considered here, with the wave equation written as a first order system, reads
\begin{align}
	u_{t}(x,t) + v_{x}(x,t) &= 0 \nonumber \\
	v_{t}(x,t) + u_{x}(x,t) &= 0 \nonumber
\end{align}
on $x \in [0,1]$ with periodic boundary conditions and
\begin{equation}
	u(x,0) = \exp\left( -\frac{1}{2} \left( \frac{x - 0.5}{0.1} \right)^{2} \right), \quad v(x,0) = 0\nonumber
\end{equation}
for $0 \leq t \leq T$. For the spatial derivatives, centered differences of $4^{\rm th}$ order with $128$ points are used on the fine level and of $2^{\rm nd}$ order with $64$ points on the coarse.
Both SDC and MLSDC perform $40$ timesteps of length $\Delta t = 0.025$ to integrate up to $T=1.0$ and iterations on each step are performed until $\left\| \rvec^k \right\|_{\infty} \leq 5 \times 10^{-8}$.
The average number of fine level sweeps over all steps for SDC and MLSDC is shown in Table~\ref{tab:wave_eq} for three different values of $M$.
In all cases, MLSDC leads to savings in terms of required fine level sweeps.
We note that for a fine level spatial resolution of only $64$ points, using spatial coarsening has a significant negative effect on the performance of MLSDC (not documented here): This suggests that for a problem which is spatially under-resolved on the finest level, further coarsening the spatial resolution within MLSDC might hurt performance, see also~\S\ref{subsec:shear_layer}.
\begin{table}[th]
\centering
\begin{tabular}{|c|c|c|} \hline
$M$ & SDC & MLSDC(1,2)  \\ \hline
$3$ & 18.5 & 11.1 \\
$5$ & 17.6 & 10.6 \\
$7$ & 14.3 & 8.2 \\ \hline
\end{tabular}
\caption{Average number of fine level sweeps over all time-steps of SDC and MLSDC for the wave equation example to reach a residual of $\left\| \rvec^{k} \right\|_{\infty} \leq 5 \times 10^{-8}$. The numbers in parentheses after MLSDC indicate the used coarsening strategies, see~\S\ref{sec:mlsdc_spatial_coarsening}.}\label{tab:wave_eq}
\end{table}

\subsection{1D viscous Burgers' equation}\label{subsec:visc_burg}
In this section we investigate the effect of coarsening in MLSDC
by considering the  nonlinear viscous Burgers' equation
\begin{align}
	\pdesol{u}_{t} + \pdesol{u}\cdot\pdesol{u}_{x} &= \nu \pdesol{u}_{xx}, \ x \in [-1,1], \ t \in [0,t_{\rm end}] \nonumber \\
		\pdesol{u}(x,0) &= u^{0}(x) \label{eq:pde} \\
		\pdesol{u}(-1,t) &= \pdesol{u}(1,t), \nonumber
\end{align}
with $\nu > 0$ and initial condition
\begin{equation}
	u^{0}(x) = \exp\left( -\frac{x^2}{\sigma^2} \right), \quad \sigma = 0.1\nonumber
\end{equation}
corresponding to a Gaussian peak strongly localized around $x=0$.
We denote the evaluation of the continuous function $\pdesol{u}$ on a given spatial mesh with points $(x_{i})_{i=1,\ldots,N}$ with a subscript $N$, so that
\begin{equation}
	\pdesol{u}_{N}(t) := \left( \pdesol{u}(x_{i}, t) \right)_{i=1,\ldots,N} \in \mathbb{R}^{N}.\nonumber
\end{equation}
Discretization of~\eqref{eq:pde} in space then yields an initial value problem
\begin{align}
	\label{eq:ivp}
	\odesol{u}_{t}(t) &= f_{N}(\odesol{u}(t)), \quad \odesol{u}(t) \in \mathbb{R}^{N}, \quad t \in [0,t_{\rm end}] \nonumber \\
	\odesol{u}(0) &= \pdesol{u}^{0}_{N}
\end{align}
with solution $u$. Finally, we denote by $\disc{u}_{N, M, \Delta t,k} \in \mathbb{R}^{N}$ the result of solving~\eqref{eq:ivp} with $k$ iterations of MLSDC using a timestep of $\Delta t$, $M$ substeps (or $M+1$ Lobatto collocation nodes), and an $N$-point spatial mesh on the finest level over one time step.

Two runs are performed here, solving~\eqref{eq:pde} with $\nu=1.0$ and $\nu=0.1$ with a single MLSDC timestep $t_{\rm end} = \Delta t = 0.01$.
MLSDC with two levels with 7 Gauss-Lobatto collocation points is used with a spatial mesh of $N = 256$ points on the fine level, and $N=128$ on the coarse level (Strategy 1).
The advective term is discretized using a $5^{\rm th}$-order WENO finite difference method~\cite{JiangShu1996} on the fine level and a simple $1^{\rm st}$-order upwind scheme on the coarse level.
For the Laplacian, a $4^{\rm th}$-order compact stencil is used on the fine level and a $2^{\rm nd}$-order stencil is used on the coarse level (Strategy 2).
The advective term is treated explicitly while the diffusion term is treated implicitly.
The resulting linear system is solved using a linear multigrid solver with a tolerance of $5\times10^{-14}$ on the fine level but solved only approximately using a single V-cycle on the coarse level (Strategy 3).
A fixed number of $K=80$ MLSDC iterations is performed here without setting a tolerance for the MLSDC residual.

In order to assess the different error components, a reference PDE solution $\pdesol{u}_{N}(\Delta t)$ is computed with a single-level SDC scheme on a mesh with $N=1,024$ points using $M+1=9$ and $\Delta t = 10^{-4}$.
An ODE solution $\odesol{u}(\Delta t)$ is computed by running single-level SDC using $M+1=9$, $\Delta t = 10^{-4}$ and the same spatial discretization as on the fine level of the MLSDC run.
Finally, the collocation solution $\odesol{u}^{\rm coll}(\Delta t)$ is computed by performing $100$ iterations of single-level SDC with $M+1=7$ and again the same spatial discretization as the MLSDC fine level.
Reference ODE and collocation solutions are computed for the coarse level using the same parameters and the MLSDC coarse level spatial discretization.

\subsubsection{Error components in MLSDC}\label{subsubsec:error_comp}
The relative error of the fully discrete MLSDC solution to the analytical solution $\pdesol{u}$ of the PDE~\eqref{eq:pde} after a single timestep of length $\Delta t$ is given by
\begin{equation}
	\label{eq:pde_error}
	\varepsilon^{\rm PDE} := \frac{ \left\| \pdesol{u}_{N}(\Delta t) - \disc{u}_{N, M, \Delta t,k} \right\| }{ \left\| \pdesol{u}_{N}(\Delta t) \right\| },
\end{equation}
where $\left\| \cdot \right\|$ denotes some norm on $\mathbb{R}^{N}$.
All errors are hereafter reported using the maximum norm $\left\| \cdot \right\|_{\infty}$.
The error $\varepsilon^{\rm PDE}$ includes contributions from three sources
\begin{align}
	\epsilon_{N} &\coloneq \frac{ \left\| \pdesol{u}_{N}(\Delta t) - \odesol{u}(\Delta t) \right\|}{\left\| \pdesol{u}_{N}(\Delta t) \right\|} \approx \text{(i) -- relative spatial error},\nonumber \\
	\epsilon_{\Delta t} &\coloneq \frac{ \left\| \odesol{u}(\Delta t) - \odesol{u}^{\rm coll}(\Delta t) \right\| }{ \left\| \pdesol{u}_{N}(\Delta t) \right\| } \approx \text{(ii) -- relative temporal error},\nonumber \\
	\varepsilon^{\rm coll} &\coloneq \frac{ \left\| \odesol{u}^{\rm coll}(\Delta t) - \disc{u}_{N,M,\Delta t,k} \right\| }{ \left\| \pdesol{u}_{N}(\Delta t) \right\|} \approx \text{(iii) -- iteration error},\nonumber
\end{align}
with $u^{\rm coll}$ denoting the exact solution of the collocation equation~\eqref{eq:compact}.
Here, (i) is the spatial discretization error; (ii) is the temporal discretization error, which is the error from replacing the analytical Picard formulation~\eqref{eq:picard} with the discrete collocation problem~\eqref{eq:compact}; and (iii) is the error from solving the collocation equation approximately using the MLSDC iteration.
The PDE error~\eqref{eq:pde_error} can be estimated using the triangle inequality according to
\begin{equation}
	\varepsilon^{\rm PDE} \leq \epsilon^{N} + \epsilon_{\Delta t} + \varepsilon^{\rm coll}.\nonumber
\end{equation}

In addition to the PDE error, we define the error between the MLSDC solution and the analytical solution of the semi-discrete ODE~\eqref{eq:ivp} as
\begin{equation}
	\label{eq:temp_error}
	\varepsilon^{\rm ODE} := \frac{ \left\| \odesol{u}(\Delta t) - \disc{u}_{N,M,\Delta t,k} \right\| }{ \left\| \pdesol{u}_{N}(\Delta t) \right\| }
		\leq \varepsilon_{\Delta t} + \varepsilon^{\rm coll}.
\end{equation}
Note that $\varepsilon^{\rm ODE}$ contains contributions from (ii) and (iii), and once the MLSDC iteration has converged, error~\eqref{eq:temp_error} reduces to the error arising from replacing the exact Picard integral~\eqref{eq:picard} by the collocation formula~\eqref{eq:compact}.

The three different error components of MLSDC, $\varepsilon^{\rm PDE}$, $\varepsilon^{\rm ODE}$ and $\varepsilon^{\rm coll}$ are expected to saturate at different levels as $k \to \infty$ according to
\begin{align}
	\varepsilon^{\rm PDE} &\to \max\{ \epsilon_{N}, \epsilon_{\Delta t} \},\nonumber  \\
	\varepsilon^{\rm ODE} &\to \epsilon_{\Delta t}, \text{ and }\nonumber  \\
	\varepsilon^{\rm coll} &\to 0.\nonumber
\end{align}
The crucial point here is that due to the presence of the FAS correction included in MLSDC, we expect $\varepsilon^{\rm PDE}$, $\varepsilon^{\rm ODE}$ and $\varepsilon^{\rm coll}$ on \emph{all} levels to saturate at values of $\epsilon_{N}$ and $\epsilon_{\Delta t}$ determined by the discretization used on the \emph{finest} level.
That is, the FAS correction should allow MLSDC to represent the solution on all coarse levels to the same accuracy as on the finest level.
This is verified in \S\ref{subsubsec:all_level_conv}.

\subsubsection{Convergence of MLSDC on all levels}\label{subsubsec:all_level_conv}
\begin{figure}[th]
	\centering
	\subfloat[Errors on fine level for $\nu = 0.1$.\label{fig:burgers_a}]{\includegraphics[width=0.48\textwidth,type=pdf,ext=.pdf,read=.pdf]{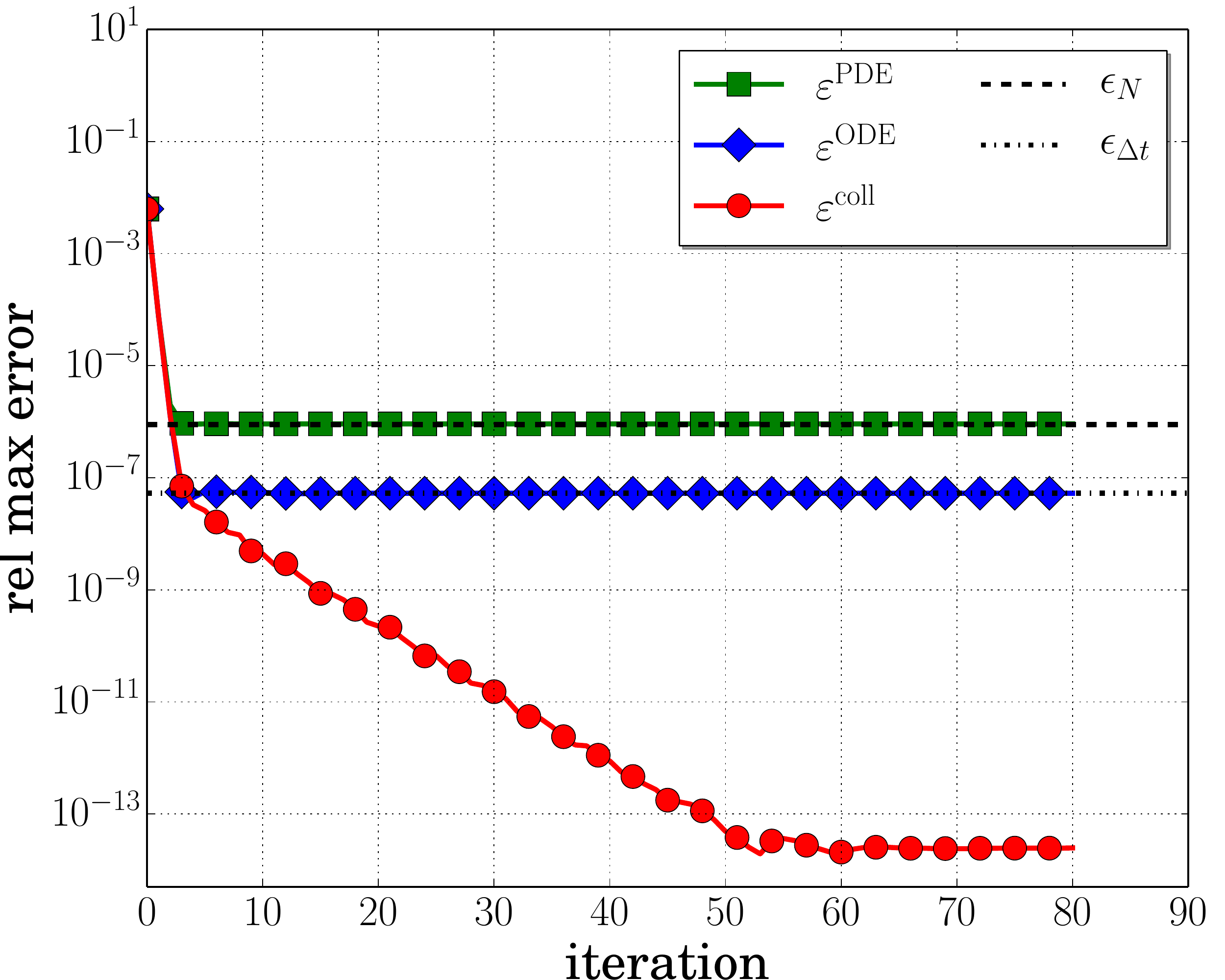}}
	\subfloat[Errors on coarse level for $\nu = 0.1$.\label{fig:burgers_b}]{\includegraphics[width=0.48\textwidth,type=pdf,ext=.pdf,read=.pdf]{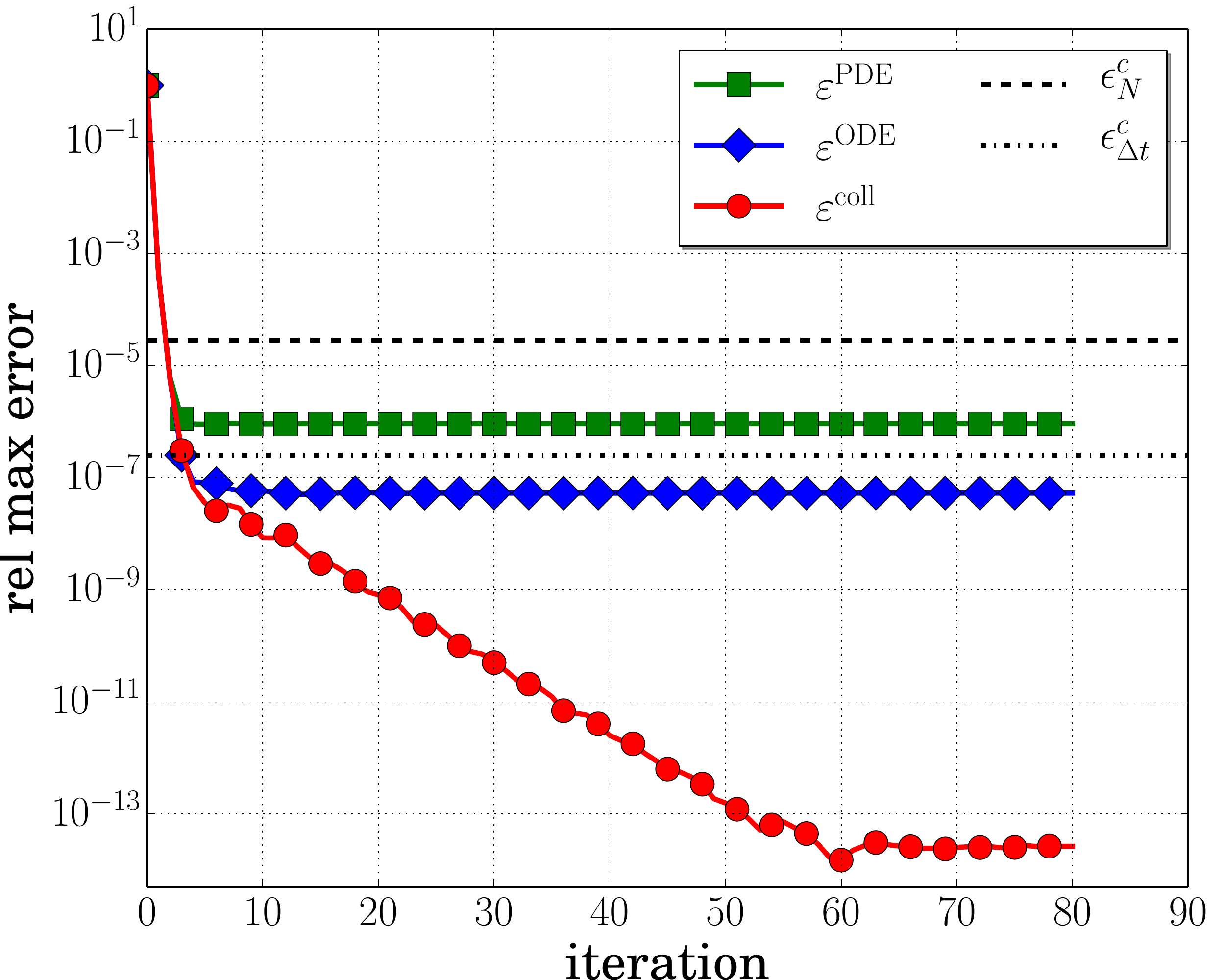}}\newline \centering
	\subfloat[Errors on fine level for $\nu = 1.0$\label{fig:burgers_c}]{\includegraphics[width=0.48\textwidth,type=pdf,ext=.pdf,read=.pdf]{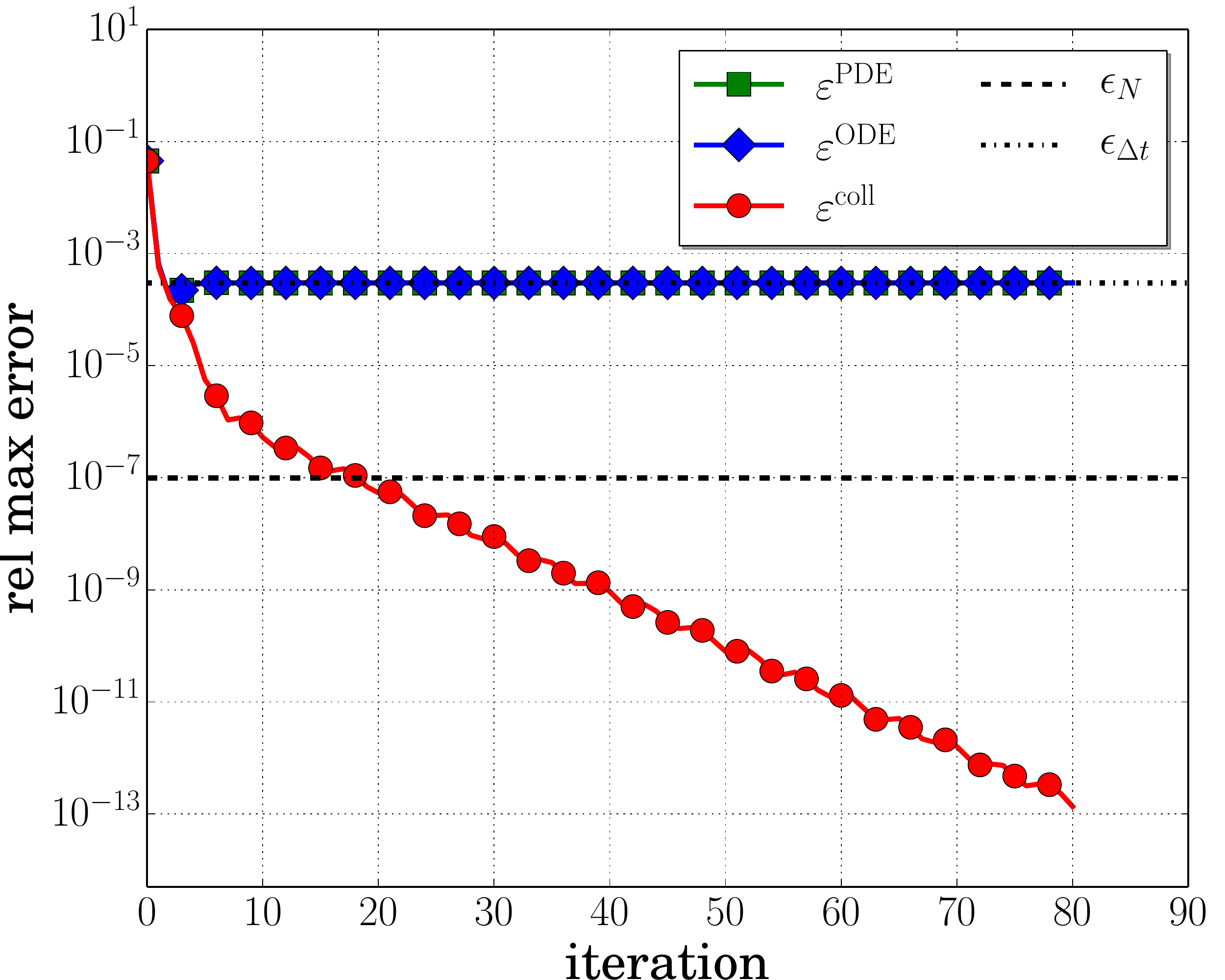}}
	\subfloat[Errors on coarse level for $\nu=1.0$\label{fig:burgers_d}]{\includegraphics[width=0.48\textwidth,type=pdf,ext=.pdf,read=.pdf]{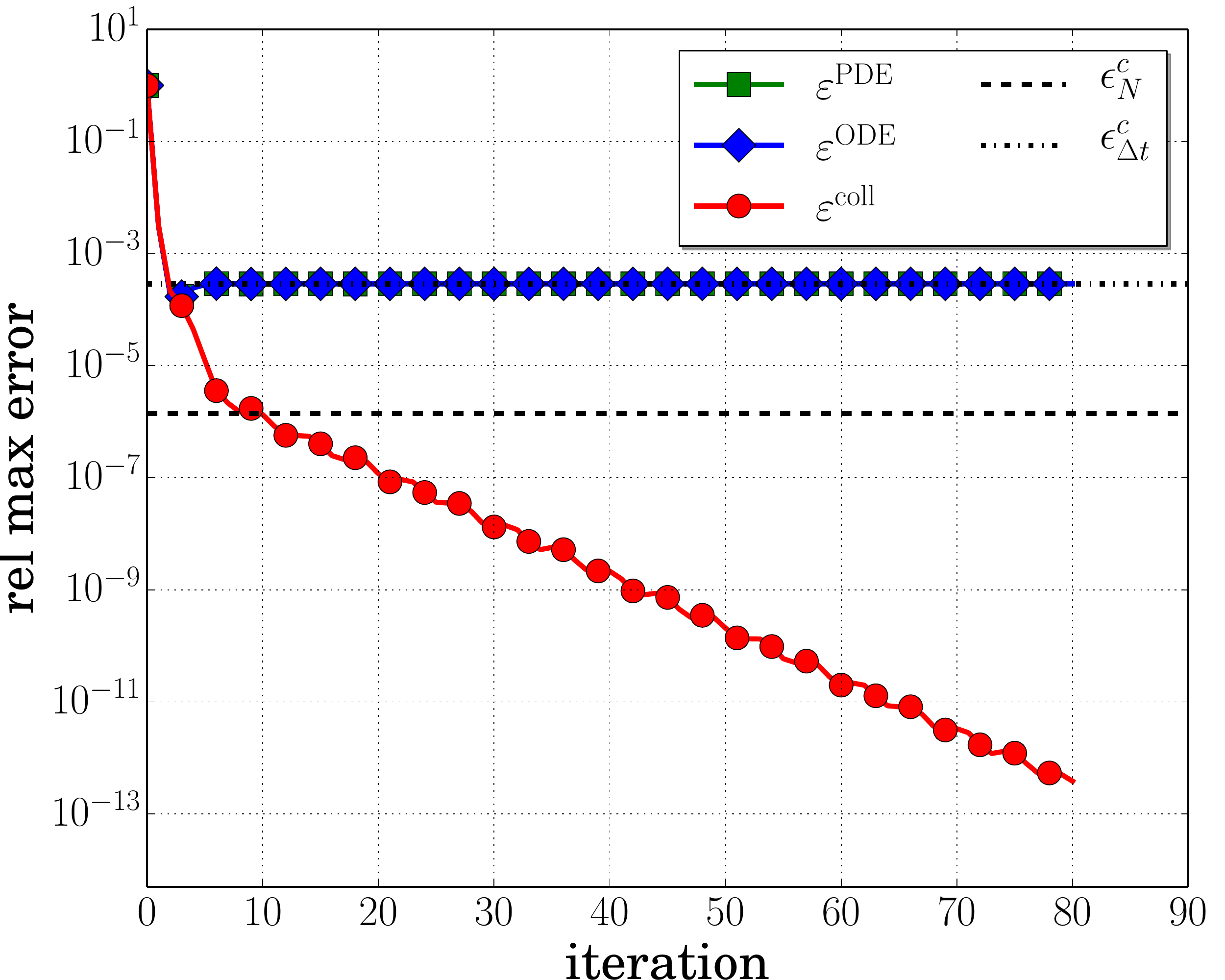}}
	\caption{Errors on fine and coarse level of MLSDC vs. iteration count.
The dashed line indicates the spatial error $\epsilon_{N}$ while the dot-dashed line indicates the temporal error $\epsilon_{\Delta t}$.
The red circles indicate the difference $\varepsilon^{\rm coll}$ between MLSDC and the collocation solution,
the blue diamonds indicate the difference $\varepsilon^{\rm ODE}$ between MLSDC and the ODE solution and the green squares indicate the difference $\varepsilon^{\rm PDE}$ between MLSDC and the PDE solution.
In (c) and (d), $\varepsilon^{\rm ODE}$ is nearly identical to $\varepsilon^{\rm PDE}$.
Note how the FAS correction in MLSDC allows the coarse level to attain the same accuracy as the fine level solution:
the saturation limits on the fine and coarse mesh are identical.}
	\label{fig:burgers}
\end{figure}

Figure~\ref{fig:burgers} shows the three error components $\varepsilon^{\rm PDE}$ (green squares), $\varepsilon^{\rm ODE}$ (blue diamonds) and $\varepsilon^{\rm coll}$ (red circles) for $\nu=0.1$ (upper) and $\nu = 1.0$ (lower) plotted against the iteration number $k$.
The errors on the fine level are shown on the left in Figures~\ref{fig:burgers_a} and~\ref{fig:burgers_c}, while errors on the coarse mesh are shown on the right.
Furthermore, the estimated spatial discretization error $\epsilon_{N}$ (dashed) and temporal discretization error $\epsilon_{\Delta t}$ (dash-dotted) are indicated by black lines.

For $\nu=0.1$, we note that the PDE error $\varepsilon^{\rm PDE}$ on the fine level (Figures~\ref{fig:burgers_a} and~\ref{fig:burgers_c}) saturates -- as expected -- at a level determined by the spatial discretization error $\epsilon_{N}$; and the ODE error $\varepsilon^{\rm ODE}$ saturates at the level of the temporal discretization error $\epsilon_{\Delta t}$.  The collocation error $\varepsilon^{\rm coll}$ saturates at near machine accuracy.
Increasing the viscosity to $\nu=1.0$, the spatial error remains at about $10^{-7}$ on the fine level but the time discretization error significantly increases compared to $\nu=0.1$.
Thus in Figure~\ref{fig:burgers_c}, both the PDE and the ODE error saturate at the value indicated by $\epsilon_{\Delta t}$.
Once again, the collocation error goes down to machine accuracy, although the rate of convergence is somewhat slower compared to $\nu=0.1$.

On the coarse level (Figures~\ref{fig:burgers_b} and~\ref{fig:burgers_d}), the estimated spatial error $\epsilon_{N}$ is noticeably higher because the values of $N$ are smaller and the order of the spatial discretization is lower.
However, as expected, the coarse level error of MLSDC saturates at values determined by the accuracy of the \emph{finest} level.  The saturation of $\varepsilon^{\rm PDE}$ and $\varepsilon^{\rm ODE}$ are identical in the left and right figures, despite the difference in $\epsilon_{N}$ and $\epsilon_{\Delta t}$.
This demonstrates that the FAS correction in MLSDC allows the solutions on coarse levels to obtain the accuracy of the finest level as long as sufficiently many iterations are performed.

\subsubsection{Required iterations}
Table~\ref{tab:burger_it} shows the number of fine level sweeps required by SDC and MLSDC to reduce the infinity norm of the residual $\rvec^k$, see ~\eqref{eq:residual}, below $10^{-5}$.
For both setups, $\nu=0.1$ as well as $\nu=1.0$, MLSDC reduces the number of required fine sweeps compared to single-level SDC.
In turn, however, MLSDC adds some overhead from coarse level sweeps.
If these are cheap enough, the reduced iteration number will result in reduced computing time, cf. \S\ref{subsec:visc_burger3d}.
\begin{table}[t]
\centering
\begin{tabular}{|c|c|c|c|}\hline
\multicolumn{2}{|c|}{$\nu=0.1$} & \multicolumn{2}{|c|}{$\nu=1.0$} \\ \hline
Method   & \# Fine sweeps & Method & \# Fine sweeps \\ \hline
SDC        & 4 & SDC         & 12 \\
MLSDC   & 3 & MLSDC  &  7 \\ \hline
\end{tabular}
\caption{Number of fine level sweeps required to reach a residual of $\left\| \rvec^{k} \right\|_{\infty} \leq 10^{-5}$ for SDC and multi-level SDC for Burgers' equation with $\nu=0.1$ and $\nu=1.0$.}\label{tab:burger_it}
\end{table}

\subsubsection{Stopping criteria} Note that the overall PDE error of the solution is not reduced further by additional iterations once $\varepsilon^{\rm coll} \leq \max\{ \epsilon_{N}, \epsilon_{\Delta t} \}$.
In Figures~\ref{fig:burgers_a}--\ref{fig:burgers_d}, this corresponds to the point where the line with red circles (iteration error) drops below the dot-dashed line (indicating $\varepsilon_{\Delta t}$) or dashed line (indicating $\varepsilon_{N}$).
The MLSDC solution, however, continues to converge to the collocation solution.
In a scenario where the PDE error is the main criterion for the quality of a solution, iterating beyond $\varepsilon^{\rm PDE}$ no longer improves the solution.
This suggests adaptively setting the tolerance for the residual of the MLSDC iteration in accordance with error estimators for $\epsilon_{\rm N}$ and $\epsilon_{\Delta t}$ to avoid unnecessary further iterations.

\subsection{Shear layer instability}\label{subsec:shear_layer}
In this example, we study the behavior of MLSDC in the case where the exact solution is not
well resolved.
We consider a shear layer instability in a 2D doubly periodic domain governed by the vorticity-velocity formulation of the 2D Navier-Stokes equations given by
\begin{align}
	\pdesol{\omega}_t + \pdesol{u}\cdot\nabla\pdesol{\omega} = \nu\nabla^{2}\pdesol{\omega}\nonumber
\end{align}
with velocity $\pdesol{u}\in\mathbb{R}^2\times[0,\infty)$, vorticity $\pdesol{\omega} = \nabla\times\pdesol{u}\in\mathbb{R}^\times[0,\infty)$ and viscosity $\nu\in\mathbb{R}^{+}$. We consider the spatial domain $[0,1]^2$ with periodic boundary conditions in all directions and the initial conditions
\begin{align}
	\pdesol{u}_1^0(x,y) &= -1.0 + \tanh(\rho(0.5 - y)) + \tanh(\rho(y - 0.25))\nonumber\\
	\pdesol{u}_2^0(x,y) &= -\delta \sin(2\pi(x + 0.25)).\nonumber
\end{align}
These initial conditions correspond to two horizontal shear layers, of ``thickness'' $\rho = 50$, at $y = 0.75$ and $y = 0.25$, with a disturbance of magnitude $\delta= 0.05$ in the vertical velocity $\pdesol{u}_2$.
As in \S\ref{subsec:visc_burg}, the system is split into implicit/explicit parts according to
\begin{align}
	\pdesol{\omega}_t &= f^E(\pdesol{\omega}) + f^I(\pdesol{\omega})\nonumber
	% \pdesol{\omega}(0) &= \pdesol{\omega}^0_N,
\end{align}
% where $\pdesol{\omega}^0_N$ is the pointwise evaluation of the continuous initial value on a $N\times N$ mesh, i.e.~the 2D analogue of~\eqref{eq:mesh_eval}, and
where
\begin{align}
	f^E(\pdesol{\omega}) &= -\pdesol{u}\cdot\nabla\pdesol{\omega}\nonumber\\
	f^I(\pdesol{\omega}) &= \nu\nabla^{2}\pdesol{\omega}.\nonumber
\end{align}
While the implicit term $f^I$ is discretized and solved as before, we apply a streamfunction approach for the explicit term $f^E$: for periodic boundary conditions, we can assume $\pdesol{u} = \nabla\times\pdesol{\psi}$ for a solenoidal streamfunction $\pdesol{\psi}$. Thus,
\begin{align}
	\pdesol{\omega} = \nabla\times(\nabla\times\pdesol{\psi}) = -\nabla^{2}\pdesol{\psi}.\nonumber
\end{align}
We refer to~\cite{chorin_mathematical_1990} for more details.
To compute $f^E_{p,N}(\odesol{\omega})$ with order-$p$ operators on an $N\times N$ mesh, we therefore solve the Poisson problem
\begin{align}
	-\nabla^{2}\pdesol{\psi} = \pdesol{\omega}\nonumber
\end{align}
for $\pdesol{\psi}$ using the linear multigrid method described previously, calculate the discretized version of $\pdesol{u} = \nabla\times\pdesol{\psi}$ and finally compute the discretization of $\pdesol{u}\cdot\nabla\pdesol{\omega}$, both with order-$p$ operators.

Two levels with $M+1=9$ collocation nodes are used with a $128 \times 128$ point spatial mesh and a fourth order stencil on the fine level.
Different combinations of coarsening are tested (the numbers in parentheses correspond to the strategies as listed in \S\ref{sec:mlsdc_spatial_coarsening}):
\begin{enumerate}
	\item MLSDC(1,2) uses a coarsened $64 \times 64$ point mesh on the coarse level and second-order stencils.
	\item MLSDC(1,2,3(1)) as MLSDC(1,2) but also solves the implicit linear systems in the coarse SDC sweep only approximately with a single V-cycle.
	\item MLSDC(1,2,3(2)) as MLSDC(1,2,3(1)) but with two V-cycles.
	\item MLSDC(2,3(1)) uses also a $128 \times 128$ point mesh on the coarse level, but second-order stencils and approximate linear solves using a single V-cycle.
\end{enumerate}
The simulation computes $256$ timesteps of MLSDC up to a final time $t=1.0$.
As reference, a classical SDC solution is computed using $1024$ timesteps with $M+1=13$ collocation nodes and the fine level spatial discretization.
Both SDC and MLSDC iterate until the residual satisfies $\left\| \rvec^k \right\|_{\infty} \leq 10^{-12}$.

\subsubsection{Vorticity field on all levels}
\begin{figure}[t!]
	\centering
%	\subfloat[SDC on a $128\times 128$ grid with $M=13$, $p=4$\label{fig:sli_field_a}]{\includegraphics[width=0.48\textwidth,type=pdf,ext=.pdf,read=.pdf]{images/slsdc_NV16384_N1024_M13_nl1_qt1_te1_dt0.00097656_np6_tc0_sc0_la0.0001_level01_vfield}}\hfill
	\subfloat[MLSDC, fine level: $128\times 128$, $p=4$\label{fig:sli_field_b}]{\includegraphics[width=0.48\textwidth,type=pdf,ext=.pdf,read=.pdf]{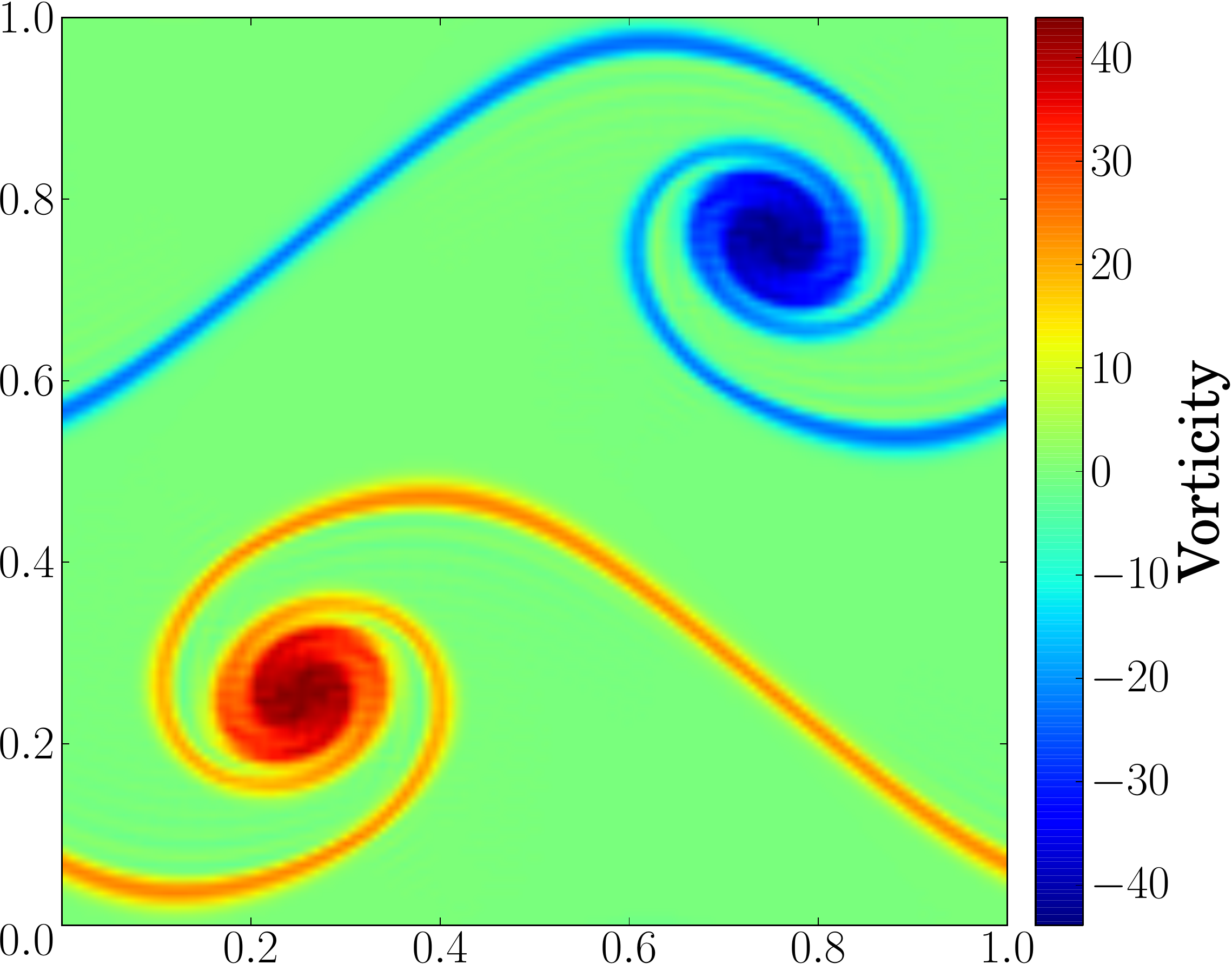}}
	\subfloat[MLSDC, coarse level: $64\times 64$, $p=2$\label{fig:sli_field_c}]{\includegraphics[width=0.48\textwidth,type=pdf,ext=.pdf,read=.pdf]{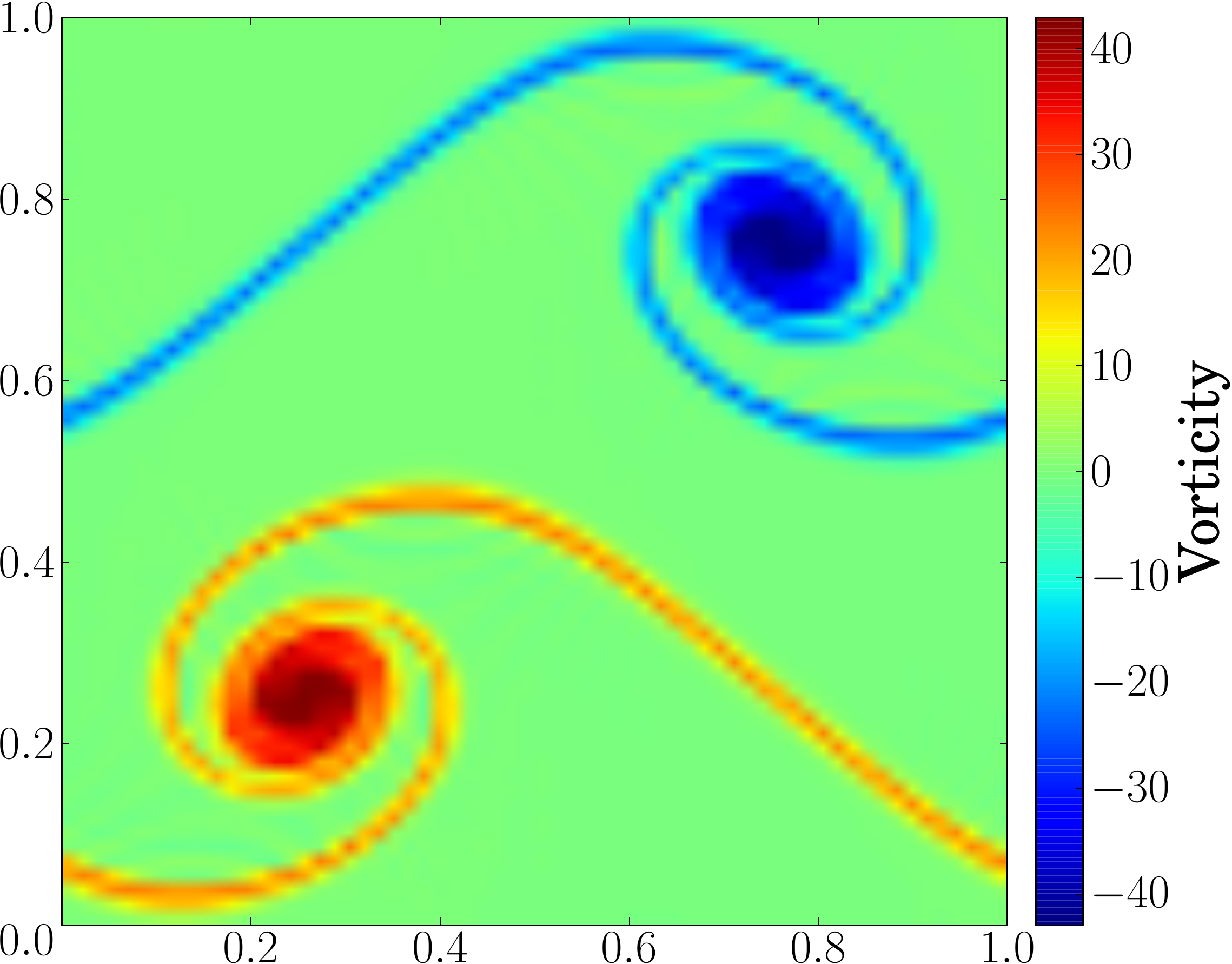}}
%	\subfloat[MLSDC, coarse level: $32\times 32$, $M=3$, $p=2$\label{fig:sli_field_d}]{\includegraphics[width=0.48\textwidth,type=pdf,ext=.pdf,read=.pdf]{images/mlsdc_NV16384_N512_M9_nl3_qt1_te2_dt0.0039062_np6_tc0_sc3_la0.0001_level03_vfield}}
	\caption{Vorticity of the solution of the shear layer instability at $t=1.0$ on the fine level (left) and coarse level (right) using MLSDC(1,2,3(1)).}
\label{fig:sli_field}
\end{figure}

Figure~\ref{fig:sli_field} shows the vorticity field at the end of the simulation on the fine and the coarse level.
The relative maximum error $\varepsilon^{\rm ODE}$ at time $t=1$ is approximately $10^{-12}$ (which corresponds to the spatial and temporal residual thresholds that were used for all runs in this example).
We note that simply running SDC with the coarse level spatial discretization from MLSDC(1,2) gives completely unsatisfactory results (not shown):  spurious vortices exist in addition to the two correct vortices and strong spurious oscillations are present in the vorticity field.
In contrast, the coarse level solution from MLSDC shown in Figure~\ref{fig:sli_field_c} looks reasonable, again because of the FAS correction.

\subsubsection{Required iterations}
Table~\ref{tab:vortex_it} shows the average number of fine level sweeps over all timesteps required by SDC and MLSDC to converge.
The configurations MLSDC(1,2), MLSDC(1,2,3(1)) and MLSDC(1,2,3(2)) do not reduce the number of sweeps, but instead lead to a small increase.
Avoiding a coarsened spatial mesh in MLSDC(2,3(1)), however, saves a small amount of fine sweeps compared to SDC.
Note that here, in contrast to the example presented in \S\ref{subsec:visc_burger3d}, Strategy 1 has a significant negative impact on the performance of MLSDC.
This illustrates that coarsening in MLSDC cannot be used in the same way  for every problem: a careful adaption of the employed strategies to the problem at hand is necessary.
\begin{table}[th]
\centering
\begin{tabular}{|c|c|} \hline
Method & \# Fine sweeps on average \\ \hline
SDC        & 6.46 \\
MLSDC(1,2) & 6.64 \\
MLSDC(1,2,3(1)) & 6.62 \\
MLSDC(1,2,3(2)) & 6.64 \\
MLSDC(2,3(1)) & 5.26 \\ \hline
\end{tabular}
\caption{Average number of fine level sweeps required to converge for SDC and MLSDC for the shear layer instability. The numbers indicate the different coarsening strategies.}\label{tab:vortex_it}
\end{table}

\subsection{Three-dimensional viscous Burgers' equation}\label{subsec:visc_burger3d}
To demonstrate that MLSDC can not only reduce iterations but also runtime, we consider viscous Burgers' equation in three dimensions
\begin{equation}
	\pdesol{u}_{t}(\textbf{x},t) + \pdesol{u}(\textbf{x},t) \cdot \nabla \pdesol{u}(\textbf{x},t)= \nu \nabla^2 \pdesol{u}(\textbf{x},t), \quad \textbf{x} \in [0,1]^{3}, \quad 0 \leq t \leq 1\nonumber
\end{equation}
with $\textbf{x} = (x,y,z)$, initial value
\begin{equation}
	\pdesol{u}(\textbf{x}, t) = \exp\left( -\frac{(x-0.5)^2 + (y-0.5)^2 + (z-0.5)^2}{\sigma^2} \right), \quad \sigma=0.1,\nonumber
\end{equation}
homogeneous Dirichlet boundary condition and diffusion coefficients $\nu=0.1$ and $\nu=1.0$.
The problem is solved using a {\sc fortran} implementation of MLSDC combined with a {\sc c} implementation of a parallel multigrid solver (PMG) in space~\cite{bolten:2014}.
A single timestep of length $\Delta t = 0.01$ is performed with MLSDC, corresponding to CFL numbers from the diffusive term on the fine level, that is
\begin{equation}
	C_{\rm diff} := \frac{\nu \Delta t}{\Delta x^{2}},\nonumber
\end{equation}
of about $C_{\rm diff}=66$ (for $\nu=0.1$) and $C_{\rm diff} = 656$ (for $\nu=1.0$).
The diffusion term is integrated implicitly using PMG to solve the corresponding linear system and the advection term
is treated explicitly.
Simulations are run on $512$ cores on the IBM BlueGene/Q JUQUEEN at the J\"ulich Supercomputing Centre.

MLSDC is run with $M+1=3$, $M+1=5$ and $M+1=7$ Gauss-Lobatto nodes with a tolerance for the residual of $10^{-5}$.
Two MLSDC levels are used with  all three types of coarsening applied:
\begin{enumerate}
\item The fine level uses a $255^{3}$ point mesh and the coarse level $127^{3}$.
\item A $4^{\rm th}$-order compact difference stencil for the Laplacian and a $5^{\rm th}$-order WENO~\cite{JiangShu1996} for the advection term are used on the fine level; a $2^{\rm nd}$-order stencil for the Laplacian and a $1^{\rm st}$-order upwind scheme for advection on the coarse.
\item The accuracy of the implicit solve on the coarse level is varied by fixing the number of V-cycles of PMG on this level.
\end{enumerate}
Three runs are performed, each with a different number of V-cycles on the coarse level.
In the first run, the coarse level linear systems are solved to full accuracy, whereas
the second and third runs use one and two V-cycles of PMG on the coarse level, respectively, instead of solving to full accuracy.  These cases are
referred to as MLSDC(1,2), MLSDC(1,2,3(1)), and MLSDC(1,2,3(2)).
On the fine level, implicit systems are always solved to full accuracy (the PMG multigrid iteration reaches a tolerance of reach a tolerance of $10^{-12}$ or stalls).
%Both SDC and MLSDC iterations are terminated once $\left\| \rvec^{k} \right\|_{\infty} \leq 10^{-5}$.

\paragraph{Required iterations and runtimes.}

Table~\ref{tab:burger3d_eq} shows both the required fine level sweeps for SDC and MLSDC as well as the total runtimes in seconds for $\nu = 0.1$ and $\nu = 1.0$ for three different values of $M$.
MLSDC(1,2) and MLSDC(1,2,3(2)) in all cases manage to significantly reduce the number of fine sweeps required for convergence in comparison to single-level SDC, typically by about a factor of two.
These savings in fine level sweeps translate into runtime savings on the order of $30-40 \%$.
For $3$ and $5$ quadrature nodes, there is no negative impact in terms of additional fine sweeps by using a reduced implicit solve on the coarse level and MLSDC(1,2,3(2)) is therefore faster than MLSDC(1,2).
However, since coarse level V-cycles are very cheap due to spatial coarsening, the additional savings in runtime are small.
For $7$ quadrature nodes, using a reduced implicit solve on the coarse level in MLSDC(1,2,3(2)) comes at the price of an additional MLSDC iteration and therefore, MLSDC(1,2) is the fastest variant in this case.

Using only a single V-cycle for implicit solves on the coarse grid in MLSDC(1,2,3(1)) results in a modest to significant increase in the number of required MLSDC iterations compared to MLSDC(1,2,3(2)) in almost all cases. The only exception is the run with $3$ nodes and $\nu=0.1$.
Therefore, MLSDC(1,2,3(1)) is typically significantly slower than MLSDC(1,2) or MLSDC(1,2,3(2)) and not recommended for use in three dimensions.
For $7$ quadrature nodes, using only a single V-cycle leads to a dramatic increase in the number of required fine sweeps and MLSDC becomes much slower than single level SDC, indicating that the inaccurate coarse level has a negative impact on convergence.
\begin{table}[t]
\centering
{\centering \bf $M+1 = 3$ Gauss-Lobatto nodes\vspace{0.5em}\linebreak}
\begin{tabular}{|c|c|c|} \hline
\multicolumn{3}{|c|}{$\nu=0.1$}   \\ \hline
Method                     & F-Sweeps   & Runtime (sec)  \\ \hline
SDC                         &  9                   & 39.4                     \\
MLSDC(1,2)           &  4                   & 26.2                    \\
MLSDC(1,2,3(2))   &  4                  &  25.6                    \\
MLSDC(1,2,3(1))   &  5                  &  29.7                    \\ \hline
\end{tabular}
\begin{tabular}{|c|c|c|} \hline
\multicolumn{3}{|c|}{$\nu=1.0$}   \\ \hline
Method                    & F-Sweeps   & Runtime (sec)  \\ \hline
SDC                         &  16   & 74.1                     \\
MLSDC(1,2)           &  8     & 49.1                     \\
MLSDC(1,2,3(2))   &  8     & 47.0                    \\
MLSDC(1,2,3(1))   &  8     & 46.7                    \\ \hline
\end{tabular}\vspace{1em}
{\centering \bf $M+1 = 5$ Gauss-Lobatto nodes\vspace{0.5em}\linebreak}
\begin{tabular}{|c|c|c|} \hline
\multicolumn{3}{|c|}{$\nu=0.1$}   \\ \hline
Method                     & F-Sweeps   & Runtime (sec)  \\ \hline
SDC                         &  7                   & 59.5                    \\
MLSDC(1,2)           &  3                   & 40.8                  \\
MLSDC(1,2,3(2))   &  3                  &  39.8                  \\
MLSDC(1,2,3(1))   &  8                  &  79.7 \\ \hline
\end{tabular}
\begin{tabular}{|c|c|c|} \hline
\multicolumn{3}{|c|}{$\nu=1.0$}   \\ \hline
Method                    & F-Sweeps   & Runtime (sec)  \\ \hline
SDC                         & 18                 & 162.7                     \\
MLSDC(1,2)           &   9                 & 105.6                   \\
MLSDC(1,2,3(2))   &   9                &  101.5                   \\
MLSDC(1,2,3(1))   &  14               &  142.8 \\ \hline
\end{tabular}\vspace{1em}
{\centering \bf $M+1 = 7$ Gauss-Lobatto nodes\vspace{0.5em}\linebreak}
\begin{tabular}{|c|c|c|} \hline
\multicolumn{3}{|c|}{$\nu=0.1$}   \\ \hline
Method                     & F-Sweeps   & Runtime (sec)  \\ \hline
SDC                         &  5                  &  82.4                   \\
MLSDC(1,2)           &  2                  &  46.1                  \\
MLSDC(1,2,3(2))   &  3                  &  57.2                   \\
MLSDC(1,2,3(1))   &  11                &  147.2  \\ \hline
\end{tabular}
\begin{tabular}{|c|c|c|} \hline
\multicolumn{3}{|c|}{$\nu=1.0$}   \\ \hline
Method                    & F-Sweeps   & Runtime (sec)  \\ \hline
SDC                         & 17                 &  224.7                   \\
MLSDC(1,2)           &   8                 &  139.5                  \\
MLSDC(1,2,3(2))   &   9                 &  148.1                  \\
MLSDC(1,2,3(1))   & 44                &  560.4                 \\ \hline
\end{tabular}
\caption{Number of required fine level sweeps and resulting runtimes in seconds by SDC and MLSDC for 3D viscous Burgers' equation. The numbers in parentheses after MLSDC indicate the employed coarsening strategies, see~\S\ref{sec:mlsdc_spatial_coarsening}. Reduced implicit solves are indicated by $3(n)$ where $n$ indicates the fixed number of multigrid V-cycles. Otherwise, PMG iterates until a residual of $10^{-12}$ is reached or the iteration stalls. The tolerance for the SDC/MLSDC iteration is $10^{-5}$.}\label{tab:burger3d_eq}
\end{table}

%% section on conclusions and outlook: links to pfasst and outlook on space-time MG
\section{Discussion} \label{sec:outro}
The paper analyzes the multi-level spectral deferred correction method (MLSDC), an extension to the original single-level spectral deferred corrections (SDC) as well as ladder SDC methods.
In contrast to SDC, MLSDC performs correction sweeps in time on a hierarchy of discretization levels, similar to V-cycles in classical multigrid.
An FAS correction is used to increase the accuracy  on coarse levels.
The paper also presents a new procedure to incorporate weighting matrices arising in higher-order
compact finite difference stencils into the SDC method.
The advantage of MLSDC is that it shifts computational work from the fine level to coarse levels, thereby reducing the number of fine SDC sweeps  and, therefore, the time-to-solution.

For MLSDC to be efficient, a reduced representation of the problem on the coarse levels has to be used in order to make coarse level sweeps cheap in terms of computing time.
Three strategies are investigated numerically, namely (1) using fewer degrees of freedom,
(2) reducing the order of the discretization, and (3) reducing the accuracy of the linear solver in implicit substeps on the coarse level.
Numerical results are presented for the wave equation, viscous Burgers' equation in 1D and 3D and for the 2D Navier-Stokes equation in vorticity-velocity formulation.
It is demonstrated that because of the FAS correction, the solutions on all levels converge up to the accuracy determined by the discretization on the finest level.
More significantly, in all four examples, MLSDC can reduce the number of fine level sweeps required to converge compared to single level SDC.
For the 3D example this translates directly into significantly reduced computing times in comparison to single-level SDC.

One potential continuation of this work is to investigate reducing the accuracy of
implicit solves on the fine level in MLSDC as well.  In~\cite{SpeckEtAl2014_DDM2013}, so called
{\it inexact} spectral deferred corrections (ISDC) methods are considered, where implicit solves at
each SDC node are replaced by a small number of multigrid V-cycles. As with MLSDC, the reduced cost ofx
implicit solves are somewhat offset by an increase in the number of SDC iterations required for
convergence.  Nevertheless, numerical results in \cite{SpeckEtAl2014_DDM2013} demonstrate an overall reduction of
cost for ISDC methods versus SDC for certain test cases.  The optimal
combination of coarsening and reducing V-cycles for SDC methods using multigrid for implicit solves
appears to be problem-dependent, and an analysis of this topic is in preparation.

The MLSDC algorithm has also been applied to Adaptive Mesh Refinement (AMR) methods popular in finite-volume methods for conservative systems.
In the AMR + MLSDC algorithm, each AMR level is associated with its own MLSDC level, resulting in a hierarchy of hybrid space/time discretizations with increasing space/time resolution. When a new (high resolution) level is added to the AMR hierarchy, a new MLSDC level is created.
The resulting scheme differs from traditional sub-cycling AMR time-stepping schemes in a few notable aspects: fine level sub-cycling is achieved through increased temporal resolution of the MLSDC nodes; flux corrections across coarse/fine AMR grid boundaries are naturally incorporated into the MLSDC FAS correction; fine AMR ghost cells eventually become high-order accurate through the iterative nature of MLSDC V-cycling; and finally, the cost of implicit solves on all levels decreases with each MLSDC V-cycle as initial guesses improve.
Preliminary results suggest that the AMR+MLSDC algorithm can be successfully applied to the compressible Navier-Stokes equations with stiff chemistry for the direct numerical simulation of combustion problems.
A detailed description of the AMR+MLSDC algorithm with applications is currently in preparation.

Finally, the impact and performance of the coarsening strategies presented here are also of relevance to the parallel full approximation scheme in space and time (PFASST)~\cite{EmmettMinion2012_DDM,EmmettMinion2012,Minion2010,SpeckEtAl2012} algorithm, which is a time-parallel scheme for ODEs and PDEs.
Like MLSDC, PFASST employs a hierarchy of levels but
performs SDC sweeps on multiple time intervals concurrently with
corrections to initial conditions being communicated forward in time during the iterations.
Parallel efficiency in PFASST can be achieved  because fine SDC sweeps are done in parallel while
sweeps on the coarsest level are in essence done serially.
In the PFASST algorithm, there is a trade-off  between decreasing the cost on coarse levels
to improve parallel efficiency and retaining good accuracy on the coarse level to minimize the
number of parallel iterations required to converge.
In~\cite{EmmettMinion2012} it was shown that, for mesh-based PDE discretizations, using a spatial mesh with fewer points  on the coarse level
in conjunction with a reduced number of quadrature nodes, led to a method with significant parallel speed up.
Incorporating the additional coarsening strategies presented here for MLSDC into PFASST would
further reduce the cost of coarse levels, but it is unclear how this might translate into
an increase in the number of parallel PFASST iterations required.

%% extensive list of acknowledgements

\begin{acknowledgements}
The plots were generated with the Python Matplotlib~\cite{Hunter2007} package. The final publication is available at springerlink.com, see 
\href{http://dx.doi.org/10.1007/s10543-014-0517-x}{http://dx.doi.org/10.1007/s10543-014-0517-x}.
\end{acknowledgements}

\bibliographystyle{spmpsci}      % mathematics and physical sciences

%\bibliography{paper_refs,sdc,Pint,Pint_Self}

\begin{thebibliography}{10}
\providecommand{\url}[1]{{#1}}
\providecommand{\urlprefix}{URL }
\expandafter\ifx\csname urlstyle\endcsname\relax
  \providecommand{\doi}[1]{DOI~\discretionary{}{}{}#1}\else
  \providecommand{\doi}{DOI~\discretionary{}{}{}\begingroup
  \urlstyle{rm}\Url}\fi

\bibitem{Alam2006}
Alam, J.M., Kevlahan, N.K.R., Vasilyev, O.V.: Simultaneous spaceÐtime adaptive
  wavelet solution of nonlinear parabolic differential equations.
\newblock Journal of Computational Physics \textbf{214}(2), 829 -- 857 (2006).
%\newblock \doi{http://dx.doi.org/10.1016/j.jcp.2005.10.009}

\bibitem{ascherPetzold}
Ascher, U.M., Petzold, L.R.: Computer Methods for Ordinary Differential
  Equations and Differential-Algebraic Equations.
\newblock SIAM, Philadelphia, PA (2000)

\bibitem{bohmerHemkerStetter:1984}
B\"ohmer, K., Hemker, P., Stetter, H.J.: The defect correction approach.
\newblock In: K.~B\"ohmer, H.J. Stetter (eds.) Defect Correction Methods.
  Theory and Applications, pp. 1--32. Springer-Verlag (1984)

\bibitem{bolten:2014}
Bolten, M.: Evaluation of a multigrid solver for 3-level {Toeplitz} and
  circulant matrices on {Blue Gene/Q}.
\newblock In: K.~Binder, G.~M\"unster, M.~Kremer (eds.) NIC Symposium 2014, pp.
  345--352. John von Neumann Institute for Computing (2014).
\newblock (to appear)

\bibitem{bourlioux2003high}
Bourlioux, A., Layton, A.T., Minion, M.L.: High-order multi-implicit spectral
  deferred correction methods for problems of reactive flow.
\newblock Journal of Computational Physics \textbf{189}(2), 651--675 (2003)

\bibitem{bouzarthMinion:2011}
Bouzarth, E.L., Minion, M.L.: A multirate time integrator for regularized
  stokeslets.
\newblock Journal of Computational Physics \textbf{229}(11), 4208--4224 (2010)

\bibitem{brandt:1977}
Brandt, A.: Multi-level adaptive solutions to boundary-value problems.
\newblock Math. Comp. \textbf{31}(138), 333--390 (1977)

\bibitem{briggs}
Briggs, W.L.: A Multigrid Tutorial.
\newblock SIAM, Philadelphia, PA (1987)

\bibitem{chorin_mathematical_1990}
Chorin, A.J., Marsden, J.E.: A mathematical introduction to fluid mechanics,
  2nd edn.
\newblock Springer-Verlag (1990)

\bibitem{chow2006}
Chow, E., Falgout, R.D., Hu, J.J., Tuminaro, R.S., Yang, U.M.: A survey of
  parallelization techniques for multigrid solvers.
\newblock In: Parallel Processing for Scientific Computing, SIAM Series of
  Software, Environements and Tools. SIAM (2006)

\bibitem{ChristliebEtAl2011_CMS}
Christlieb, A., Morton, M., Ong, B., Qiu, J.M.: Semi-implicit integral deferred
  correction constructed with additive {R}unge--{K}utta methods.
\newblock Communications in Mathematical Science \textbf{9}(3), 879--902
  (2011).
%\newblock \urlprefix\url{http://dx.doi.org/10.4310/CMS.2011.v9.n3.a10}

\bibitem{ChristliebEtAl2009}
Christlieb, A., Ong, B., Qiu, J.M.: Comments on high-order integrators embedded
  within integral deferred correction methods.
\newblock Communications in Applied Mathematics and Computational Science
  \textbf{4}(1), 27--56 (2009).
%\newblock \urlprefix\url{http://dx.doi.org/10.2140/camcos.2009.4.27}

\bibitem{ChristliebEtAl2010_MoC}
Christlieb, A., Ong, B.W., Qiu, J.M.: Integral deferred correction methods
  constructed with high order {R}unge-{K}utta integrators.
\newblock Mathematics of Computation \textbf{79}, 761--783 (2010).
%\newblock \urlprefix\url{http://dx.doi.org/10.1090/S0025-5718-09-02276-5}

\bibitem{DaiEtAl2013_ESAIM}
Dai, X., Le~Bris, C., Legoll, F., Maday, Y.: Symmetric parareal algorithms for
  hamiltonian systems.
\newblock ESAIM: Mathematical Modelling and Numerical Analysis \textbf{47},
  717--742 (2013).
%\newblock \urlprefix\url{http://dx.doi.org/10.1051/m2an/2012046}

\bibitem{pereyra:1968}
Daniel, J.W., Pereyra, V., Schumaker, L.L.: Iterated deferred corrections for
  initial value problems.
\newblock Acta Cient. Venezolana \textbf{19}, 128--135 (1968)

\bibitem{dutt2000spectral}
Dutt, A., Greengard, L., Rokhlin, V.: Spectral deferred correction methods for
  ordinary differential equations.
\newblock BIT Numerical Mathematics \textbf{40}(2), 241--266 (2000)

\bibitem{EmmettMinion2012_DDM}
Emmett, M., Minion, M.L.: Efficient implementation of a multi-level parallel in
  time algorithm.
\newblock In: Proceedings of the 21st International Conference on Domain
  Decomposition Methods, Lecture Notes in Computational Science and Engineering
  (2012).
%\newblock \urlprefix\url{http://dd21.inria.fr/pdf/emmett_mini_15.pdf}.
\newblock (In press)

\bibitem{EmmettMinion2012}
Emmett, M., Minion, M.L.: Toward an efficient parallel in time method for
  partial differential equations.
\newblock Communications in Applied Mathematics and Computational Science
  \textbf{7}, 105--132 (2012).
%\newblock \urlprefix\url{http://dx.doi.org/10.2140/camcos.2012.7.105}

\bibitem{HagstromZhou2006}
Hagstrom, T., Zhou, R.: On the spectral deferred correction of splitting
  methods for initial value problems.
\newblock Communications in Applied Mathematics and Computational Science
  \textbf{1}(1), 169--205 (2006).
%\newblock \urlprefix\url{http://dx.doi.org/10.2140/camcos.2006.1.169}

\bibitem{HairerI}
Hairer, E., Norsett, S.P., Wanner, G.: Solving Ordinary Differential Equations
  I, Nonstiff Problems.
\newblock Springer-Verlag, Berlin (1987)

\bibitem{HairerII}
Hairer, E., Wanner, G.: Solving Ordinary Differential Equations {II}, Stiff and
  Differential-Algebraic Problems.
\newblock Springer-Verlag, Berlin (1991)

\bibitem{hansen2006convergence}
Hansen, A.C., Strain, J.: Convergence theory for spectral deferred correction.
\newblock Preprint  (2006)

\bibitem{HautWingate2013}
Haut, T., Wingate, B.: An asymptotic parallel-in-time method for highly
  oscillatory {PDE}s.
\newblock SIAM Journal on Scientific Computing  (2014).
\newblock In press

\bibitem{huang2006accelerating}
Huang, J., Jia, J., Minion, M.: Accelerating the convergence of spectral
  deferred correction methods.
\newblock Journal of Computational Physics \textbf{214}(2), 633--656 (2006)

\bibitem{HuangEtAl2006}
Huang, J., Jia, J., Minion, M.: Accelerating the convergence of spectral
  deferred correction methods.
\newblock Journal of Computational Physics \textbf{214}(2), 633 -- 656 (2006).
%\newblock \urlprefix\url{http://dx.doi.org/10.1016/j.jcp.2005.10.004}

\bibitem{Hunter2007}
Hunter, J.D.: Matplotlib: A 2{D} graphics environment.
\newblock Computing In Science \& Engineering \textbf{9}(3), 90--95 (2007)

\bibitem{JiangShu1996}
Jiang, G.S., Shu, C.W.: Efficient implementation of weighted {ENO} schemes.
\newblock Journal of Computational Physics \textbf{126}, 202--228 (1996)

\bibitem{Layton2009}
Layton, A.T.: On the efficiency of spectral deferred correction methods for
  time-dependent partial differential equations.
\newblock Applied Numerical Mathematics \textbf{59}(7), 1629 -- 1643 (2009).
%\newblock \urlprefix\url{http://dx.doi.org/10.1016/j.apnum.2008.11.004}

\bibitem{layton2004conservative}
Layton, A.T., Minion, M.L.: Conservative multi-implicit spectral deferred
  correction methods for reacting gas dynamics.
\newblock Journal of Computational Physics \textbf{194}(2), 697--715 (2004)

\bibitem{laytonMinion:2005}
Layton, A.T., Minion, M.L.: Implications of the choice of quadrature nodes for
  {P}icard integral deferred corrections methods for ordinary differential
  equations.
\newblock BIT Numerical Mathematics \textbf{45}, 341--373 (2005)

\bibitem{lele_compact_1992}
Lele, S.K.: Compact finite difference schemes with spectral-like resolution.
\newblock Journal of Computational Physics \textbf{103}(1), 16--42 (1992)

\bibitem{minion2003semi}
Minion, M.L.: Semi-implicit spectral deferred correction methods for ordinary
  differential equations.
\newblock Communications in Mathematical Sciences \textbf{1}(3), 471--500
  (2003)

%\newblock \urlprefix\url{http://projecteuclid.org/euclid.cms/1250880097}

\bibitem{minion2004semi}
Minion, M.L.: Semi-implicit projection methods for incompressible flow based on
  spectral deferred corrections.
\newblock Applied numerical mathematics \textbf{48}(3), 369--387 (2004)

\bibitem{Minion2010}
Minion, M.L.: A hybrid parareal spectral deferred corrections method.
\newblock Communications in Applied Mathematics and Computational Science
  \textbf{5}(2), 265--301 (2010).
%\newblock \urlprefix\url{http://dx.doi.org/10.2140/camcos.2010.5.265}

\bibitem{pereyra:1967}
Pereyra, V.: Iterated deferred corrections for nonlinear operator equations.
\newblock Numerische Mathematik \textbf{10}, 316--323 (1966)

\bibitem{pereyra:1966}
Pereyra, V.: On improving an approximate solution of a functional equation by
  deferred corrections.
\newblock Numerische Mathematik \textbf{8}, 376--391 (1966)

\bibitem{RuprechtKrause2012}
Ruprecht, D., Krause, R.: Explicit parallel-in-time integration of a linear
  acoustic-advection system.
\newblock Computers \& Fluids \textbf{59}(0), 72 -- 83 (2012).
%\newblock \urlprefix\url{http://dx.doi.org/10.1016/j.compfluid.2012.02.015}

\bibitem{south1977}
South, J.C., Brandt, A.: Application of a multi-level grid method to transonic
  flow calculations.
\newblock In: Transonic flow problems in turbomachinery, pp. 180--206.
  Hemisphere (1977)

\bibitem{SpeckEtAl2012}
Speck, R., Ruprecht, D., Krause, R., Emmett, M., Minion, M., Winkel, M.,
  Gibbon, P.: A massively space-time parallel {N}-body solver.
\newblock In: Proceedings of the International Conference on High Performance
  Computing, Networking, Storage and Analysis, SC '12, pp. 92:1--92:11. IEEE
  Computer Society Press, Los Alamitos, CA, USA (2012).
%\newblock \urlprefix\url{http://dx.doi.org/10.1109/SC.2012.6}

\bibitem{SpeckEtAl2014_DDM2013}
Speck, R., Ruprecht, D., Minion, M., Emmett, M., Krause, R.: Inexact spectral
  deferred corrections using single-cycle multigrid (2014).
%\newblock \urlprefix\url{http://arxiv.org/abs/1401.7824}.
\newblock ArXiv:1401.7824 [math.NA]

\bibitem{spotz_high-order_1996}
Spotz, W.F., Carey, G.F.: A high-order compact formulation for the {3D} poisson
  equation.
\newblock Numerical Methods for Partial Differential Equations \textbf{12}(2),
  235--243 (1996)

\bibitem{stetter:1974}
Stetter, H.J.: Economical global error estimation.
\newblock In: R.A. Willoughby (ed.) Stiff Differential Systems, pp. 245--258
  (1974)

\bibitem{trottenberg_multigrid:_2000}
Trottenberg, U., Oosterlee, C.W.: Multigrid: Basics, Parallelism and
  Adaptivity.
\newblock Academic Press (2000)

\bibitem{Weiser2013}
Weiser, M.: Faster {SDC} convergence on non-equidistant grids with {DIRK}
  sweeps (2013).
\newblock {ZIB} Report 13--30

\bibitem{ShuEtAl2007}
Xia, Y., Xu, Y., Shu, C.W.: Efficient time discretization for local
  discontinuous galerkin methods.
\newblock Discrete and Continuous Dynamical Systems -- Series B \textbf{8}(3),
  677 -- 693 (2007).
%\newblock \urlprefix\url{http://dx.doi.org/10.3934/dcdsb.2007.8.677}

\bibitem{zadunaisky:1964}
Zadunaisky, P.: A method for the estimation of errors propagated in the
  numerical solution of a system of ordinary differential equations.
\newblock In: G.~Contopoulus (ed.) The Theory of Orbits in the Solar System and
  in Stellar Systems. Proceedings of International Astronomical Union,
  Symposium 25 (1964)

\end{thebibliography}

\end{document}